\documentclass[onefignum,onetabnum]{siamart220329}

\usepackage{lipsum}
\usepackage{amsfonts}
\usepackage{graphicx}
\usepackage{epstopdf}
\usepackage{algorithmic}
\ifpdf
  \DeclareGraphicsExtensions{.eps,.pdf,.png,.jpg}
\else
  \DeclareGraphicsExtensions{.eps}
\fi

\usepackage{enumitem}
\setlist[enumerate]{leftmargin=.5in}
\setlist[itemize]{leftmargin=.5in}

\newsiamremark{remark}{Remark}
\newsiamremark{hypothesis}{Hypothesis}
\crefname{hypothesis}{Hypothesis}{Hypotheses}
\newsiamthm{claim}{Claim}

\headers{Compact Representations}{Johannes J. Brust}

\title{Useful Compact Matrices for Data-Fitting\thanks{Version of fall 2024.
\funding{This work was partially funded by the startup fund at Arizona State University.}}}

\author{Johannes J. Brust\thanks{School of Mathematical and Statistical Sciences, Arizona State University, Tempe, AZ 
  (\email{jjbrust@asu.edu}).}}

\usepackage{amsopn}
\usepackage{tikz}
\usetikzlibrary{positioning}

\usepackage{multirow}

\newcommand{\bmat}[1]{\begin{bmatrix}#1\end{bmatrix}} 

\renewcommand{\k}[1]{{#1}_k}
\newcommand{\ko}[1]{{#1}_{k+1}}
\newcommand{\kot}[1]{{#1}^T_{k+1}}
\newcommand{\kt}[1]{{#1}^T_k}

\newcommand{\bfgs}{{\small BFGS }}
\newcommand{\lbfgs}{{\small L-BFGS }}
\newcommand{\subsc}[2]{{#1}_{#2}}
\newcommand{\subsct}[2]{{#1}^T_{#2}}
\newcommand{\subsup}[3]{{#1}^{(#3)}_{#2}}
\newcommand{\subex}[3]{{#1}^{#3}_{#2}}
\newcommand{\subext}[3]{({#1}^{#3}_{#2})^T}
\renewcommand{\t}[1]{{#1}^T}
\newcommand{\ex}[2]{{#1}^{#2}}

\ifpdf
\hypersetup{
  pdftitle={Compact Representations},
  pdfauthor={Johannes J. Brust}
}
\fi

\begin{document}

\maketitle

\begin{abstract}
For minimization problems without 2nd derivative information, methods that estimate Hessian matrices can be very effective.
However, conventional techniques generate dense matrices that are prohibitive for large problems.
Limited-memory compact representations express the dense arrays in terms of a low rank representation
and have become the state-of-the-art for software implementations on large deterministic problems. We 
develop new compact representations that are parameterized by a choice of vectors and that reduce
to existing well known formulas for special choices. We demonstrate effectiveness of the 
compact representations for large eigenvalue computations, tensor factorizations and nonlinear regressions.
\end{abstract}

\begin{keywords}
compact representation, limited-memory, low-rank updates, canonical polyadic decomposition, matrix recurrence, quasi-Newton, trust-region, line-search, {\small SGD}, {\small L-BFGS}, eigendecomposition 
\end{keywords}

\begin{MSCcodes}
65F05, 65F55, 68U05, 15A23, 15A69, 90C06, 90C15, 90C30, 90C53
\end{MSCcodes}

\section{Introduction}
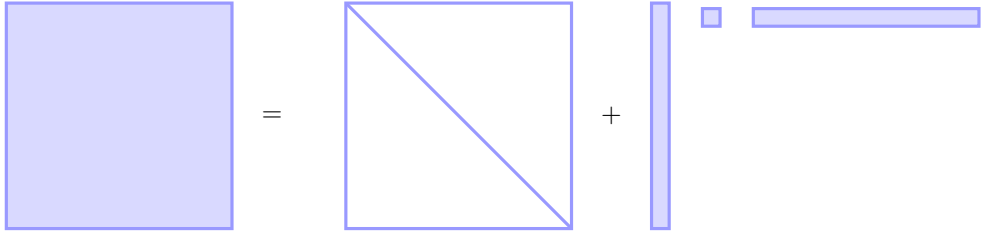
\begin{figure}
    \centering
    \begin{tikzpicture}[
    squareL/.style={rectangle, draw=blue!40, fill=blue!15, very thick, minimum height=30mm, minimum width=30mm},
    squareS/.style={rectangle, draw=blue!40, fill=blue!15, very thick, minimum size=2mm},
    rectT/.style={rectangle, draw=blue!40, fill=blue!15, very thick, minimum height=30mm, minimum width =2mm},
    rectW/.style={rectangle, draw=blue!40, fill=blue!15, very thick, minimum height=2mm, minimum width =30mm},
    init/.style={rectangle, minimum height=40.mm, minimum width =0.1mm},
    squareLE/.style={rectangle, draw=blue!40, fill=blue!0, very thick, minimum height=30mm, minimum width=30mm}
    ]
    \node[squareL]                  (main)                          {};
    \node                           (eq)        [right=2.5mm of main]     {$=$};
    \node[squareLE,path picture={ \draw[blue!40, very thick] (path picture bounding box.north west) -- (path picture bounding box.south east); }] (diag)      [right=7mm of eq]       {};
    \node                           (add)       [right=2.5mm of diag]     {$+$};                                
    \node[rectT]                    (fact1)     [right=2.5mm of add]     {};
    \node[squareS]                  (fact2)     [above right= -3.5mm and 4mm of fact1]     {};
    \node[rectW]                    (fact3)     [right=4mm of fact2]     {};    
    \end{tikzpicture}
    \caption{The low-rank form of a compact representation for a dense Hessian approximation.}
    \label{fig:compact}
\end{figure}

For large-scale data fitting one typically solves problems of the form
\begin{equation}
    \label{eq:main}
        \underset{ w \in \mathbb{R}^d }{ \text{ minimize } } f(w)
\end{equation}
where $ f: \mathbb{R}^d \rightarrow \mathbb{R} $ represents a loss, objective or penalty function and $w$ 
is a vector of parameters. Often the goal is to match a model's output as closely as possible to 
certain observed data, and therefore minimize the error between the two. Specific examples are e.g.,
tensor decompositions that minimize the distance between the factorization and data (Acar, Dunlavy and Kolda \cite{AcarDK11}),
logistic regressions for machine learning (Malouf \cite{Malouf02}), or nonlinear least-squares for model
calibration in science and engineering (Dennis, Gay and Walsh \cite{DennisGW81}), among others.
Typically, the gradient vector $ \nabla f(w) = g(w) $ is available, however higher derivatives are not.
Depending on the specific application, a variety of methods emerged as the de-facto                                             
standards. Specifically, for large deep neural networks stochastic gradient methods, like Adam (Kingma and Ba \cite{KingmaB14})
or versions of Stochastic Gradient Descent (SGD) (see e.g., Ma, Bassily and Belkin \cite{MaBB18}) are very effective,
as a consequence of the statistical nature of such problems. For nonlinear least-squares,
special methods such as NL2SOL of Dennis, Gay and Walsh \cite{dennis1981adaptive} can be effective, but for general large deterministic minimization tasks 
the limited memory BFGS (L-BFGS-B), method of Zhu, Byrd and Nocedal \cite{ZhuByrdNocedal97} is very popular.
Because it is applicable to general problems as in eq. \eqref{eq:main}, it is also the go-to-method when computing
tensor factorizations through minimization \cite{TensorTBX}. For large optimization, one typically exploits certain
structures of the problem or designs judicious methods in order to effectively handle the large dimensions.
For instance, when the problems are large and sparse and 2nd derivatives are available, 
then sparse Newton methods can converge rapidly, while maintaining efficient computations,
(see e.g., Gill, Murray, Saunders and Wright \cite{gillMSW84} or Coleman \cite{coleman84}). These methods are also 
applicable to large learning problems, such as recommender systems \cite{LuoEtAl15}. Yet, in machine learning
with very large amounts of data, the objective function typically has an additional structure. Namely,
the objective is typically a sum of individual functions that each are associated with an arbitrary data
point. Methods that randomly select subsets of the functions at every iteration have proven to be very
effective in practice. Therefore, methods that exploit the stochasticity, like variations of 
Stochastic Gradients \cite{LiuW15}, are the default techniques in this domain.
For dense problems, without stochastic properties, limited-memory quasi-Newton methods are widely used.
Even though the underlying problems can have a dense Hessian matrix, a limited memory method estimates
this structure with a low-rank representation, known as the so-called \emph{compact representation} (see Fig. \ref{fig:compact})
An important property is that the compact representation uses $ \mathcal{O}(d) $ memory to represent
the $d^2$ elements in the estimate of the Hessian matrix. Because of this decomposition, operations
like matrix vector products, linear system solves or eigenvalue computations can be performed with
complexity that is linear in the number of variables.

\subsection{Notation}
\label{sec:notation}
We use Householder notation so that lower case greek symbols are scalars, lower Roman letters are column vectors and upper
case Roman letters denote matrices: $\alpha, a, A$. The main iteration index is $k \ge 0$. We denote orthogonal 
and orthonormal matrices by $P$ and $Q$, interchangeably. $R$ is reserved for upper triangular matrices,
while $L$ and $D$ represent lower triangular and diagonals. At times we will use superscripts on a matrix,
to make its relation to another matrix explicit. For instance, $ L^{VY} $ represents a lower triangular matrix in
relation to $V$ and $S$. An underline below a matrix $Y$ means that its first column is removed $ \underline{Y} $,
while an overline represents the matrix without its first row $ \overline{Y} $.
We denote the identity matrix by $I$ with columns $ e_i $, and with dimensions that depend on the context.

\section{Unconstrained Optimization}
\label{sec:steps}
Many of the methods to generate iterates for minimizing eq. \ref{eq:main} can be described in terms of a symmetric matrix $\k{B} \in \mathbb{R}^{d \times d} $,
a vector (which is related to the gradient and may be stochastic) $\k{h}$, and a set of nonnegative scalars: $\k{\alpha}, \k{\sigma},  $ and $\k{\delta}$.
At times it is most effective to work directly with the inverse $\k{H} = \k{B}^{-1} $.
The matrix $\k{B}$ may represent the Hessian of the objective, an estimate thereof, or it could be a diagonal. With this convention,
Newton's method, quasi-Newton or gradient based methods like Adam can be described. In the stochastic setting the typical update
is of the form 
\begin{equation}
    \label{eq:stochup}
    \ko{w} = \k{w} - \k{\alpha}\k{H} \k{h}, \text{ where } \k{h} \text{ is stochastic. }
\end{equation}
For instance, by setting in eq. \eqref{eq:stochup} the inverse Hessian as the identity $\k{H} = I $, using a gradient associated with
one random data point i.e., $ \k{h} = \nabla f(\k{w};\text{``data}_i\text{''}) = g_i(\k{w}) $, and a small or decreasing $\k{\alpha}$ yields 
the Stochastic Gradient Descent (SGD) iteration.
When having nonrandom access to the objective, deterministic nonlinear unconstrained optimization approaches typically use line-searches (e.g. Zhang and Hager \cite{ZhangH04}) 
and/or trust-region strategies (e.g., Conn, Gould and Toint \cite{ConGT00a}). In this setting, typically the
vector $\k{h}$ is equal to the gradient $ \k{h} = g(\k{w}) = \k{g} $.
The steps in the two respective methods are determined in relation to the scalars: $ \k{\alpha} $ and $\k{\sigma}$:
\begin{align*}
\ko{w}                   &= \k{w} - \k{\alpha}\k{H}\k{h}    & &    & & \text{(line-search)}  \\                         
& && && \\ 
\text{ Solve } (\k{B} + \k{\sigma} I )\k{s}  &= - \k{h}     & &     & & \text{(trust-region)} \\ 
\ko{w} &= \k{w} + \k{s}                                     & &     & & 
\end{align*}
In line-search methods one determines desirable step lengths by approximately solving the one dimensional problem
$ \k{\alpha} = \min_{\alpha} f(\k{w} - \alpha \k{H} \k{h}) $. In trust-region methods one solves a sequence of 
shifted systems in order find a shift $\sigma \ge 0 $ that satisfies: $ (\k{B} + \sigma I)\k{s} = - \k{h} $ subject to 
$ \| \k{s} \| \le \k{\delta}$ and $ \k{B} + \sigma I \succeq 0 $. There are further details for practical
line-search and trust-region methods, but iterates are broadly selected to generate sufficient function reductions. 
Independent of which strategy is chosen, it is valuable to exploit structure 
in $\k{B}$ and its inverse $\k{H}$ for effective computation. Based on an initial matrix, the family of quasi-Newton methods uses
low-rank updates (typically, rank-1 or rank-2) to maintain an estimate of the Hessian or its inverse. In particular,
for two $d$-dimensional vectors and a symmetric initialization a matrix recursion specifies all remaining updates. Traditionally,
the vectors and initial matrix have been defined by
\begin{equation}
    \label{eq:sy}
    \k{s} = \ko{w} - \k{w}, \quad \k{y} = \ko{g} - \k{g}, \quad  \subsc{H}{0} \in \mathbb{R}^{d \times d} \text{ (symmetric)}.
\end{equation}
However, when the objective function is stochastic, differences between gradients such as $ \ko{g} - \k{g} $ (or $ \ko{w} - \k{w} $) may 
be noisy. Therefore, in such cases further options for how to choose $\k{y}$ have been introduced by Byrd et \emph{al.} \cite{ByrdHNS16}. It will become
clear later that the methods that we propose in this work are independent of the actual choice for $\k{s}$ and $\k{y}$. 
The recursion of the inverse \bfgs \cite{Broyden70,Fletcher70,Goldfarb70,Shanno70} matrix is
\begin{equation}
    \label{eq:recIBFGS}
    \ko{H} = \big(I - \k{\rho} \k{s} \k{y}^T \big) \k{H} \big(I - \k{\rho} \k{y} \k{s}^T \big) + \k{\rho} \k{s} \k{s}^T, \quad \k{\rho} = \frac{1}{ \k{s}^T \k{y} }
\end{equation}
Note that the matrix generated by this process is typically dense, with a pattern like the left hand side of Fig. \ref{fig:compact}.
So long $\k{s}^T\k{y} > 0$ for all iterations (and the initialization is positive definite), the sequence of matrices 
formed in eq. \eqref{eq:recIBFGS} are also all positive definite.

\subsection{Compact Representation}
\label{sec:compact}
Remarkably, by unwinding the recursion in eq. \eqref{eq:recIBFGS}, a closed matrix formula has been shown to exist
in Byrd, Nocedal and Schnabel \cite{ByrNS94}. By collecting the vectors $ \{ (\subsc{s}{i},\subsc{y}{i}) \}_{i=0}^{k-1} $
into matrices and defining a diagonal, a (strictly) lower triangular, and an upper triangular matrix one
defines the basic components of this formula
\begin{equation}
    \label{eq:SYD}
        \k{S} = \bmat{ \subsc{s}{0} & \ldots & \subsc{s}{k-1} }, \quad 
        \k{Y} = \bmat{ \subsc{y}{0} & \ldots & \subsc{y}{k-1} }, \quad 
\end{equation}
\begin{equation}
    \label{eq:LR}
    (\k{D})_{ij} = \subsct{s}{i} \subsc{y}{j} \: \: \text{for} \: i = j, \quad 
    (\k{L})_{ij} = \subsc{s}{i}^T \subsc{y}{j} \: \: \text{for} \: i > j, \quad
    (\k{R})_{ij} = \subsc{s}{i}^T \subsc{y}{j} \: \: \text{for} \: i \le j
\end{equation}
With the definitions in eqs. \eqref{eq:SYD} and \eqref{eq:LR}, a symmetric positive definite initialization
and positive $\subsct{s}{i}\subsc{y}{i} > 0 $ the compact representation of the matrix recursion \eqref{eq:recIBFGS} 
is \cite[Theorem 2.2]{ByrNS94}:
\begin{equation}
    \label{eq:compactIBFGS}
    \k{H} = \subsc{H}{0} + 
    \bmat{\k{S} & \subsc{H}{0} \k{Y} }
    \bmat{\k{R}^{-T}( \k{D} + \k{Y}^T \subsc{H}{0} \k{Y} ) \k{R}^{-1} & - \k{R}^{-T} \\ - \k{R}^{-1} & 0 }
    \bmat{\kt{S} \\ \kt{Y}\subsc{H}{0} }
\end{equation}
The main use of this formula is for the limited-memory setting where $ l \ll d $ denotes the memory parameter,
with typical values around $l=5$ (see e.g., \cite[Section 9]{BrustMPS22}). Then, instead
of storing the history of all vectors one limits this to the $l$ most recent pairs $ \{ (\subsc{s}{i}, \subsc{y}{i}) \}_{i=k-l}^{k-1} $ .
Further, typically the initialization is chosen as an adaptive multiple of the identity $ \subsup{H}{0}{k} = \k{\gamma} I   $,
with $ \k{\gamma} = \subsct{y}{k-1} \subsc{s}{k-1} / \subsct{y}{k-1} \subsc{y}{k-1} $. Limited-memory methods are frequently used for 
large-scale problems with many variables (i.e., $d$ can be large), in which the limited-memory matrices $ \k{S} \in \mathbb{R}^{d \times l} $ 
and $\k{Y} \in \mathbb{R}^{d \times l} $ are tall and very skinny:
$ \k{S} = \bmat{\subsc{s}{k-l-1}, \ldots, \subsc{s}{k-1}}  $ and $ \k{Y} = \bmat{\subsc{y}{k-l-1}, \ldots, \subsc{y}{k-1}}  $. 
The factorization pattern of a limited-memory representation of eq. \eqref{eq:compactIBFGS} corresponds to the right hand side
of Fig. \ref{fig:compact}. Besides being useful for constrained problems, the low rank representation is a significant reason why the compact representation is
implemented in state-of-the-art software packages, such as {\small KNITRO} \cite{ByrdNW06}  or {\small L-BFGS-B} \cite{ZhuByrdNocedal97} and large-scale
trust-region methods like {\small L-SR1} \cite{BrustBEM22}. Note that to obtain $ \k{B} $ from the compact representation
of $ \k{H} $ one can apply the Sherman-Morrison-Woodbury inverse to eq. \eqref{eq:compactIBFGS}.
Even though the compact representation is derived from unwinding recursive matrix updates, such as in eq. \eqref{eq:recIBFGS},
the number of known compact representations is modest. A main reason for this is that unwinding 
the recurrences involves nonlinear matrix relations and terms that may appear unintuitive
(see, for example, the inverses in the middle matrix of eq. \eqref{eq:compactIBFGS}).
Nevertheless, the compact representation of the {\small SR-1} update, as well as one of Broyden's updates
for systems of nonlinear equations have been derived in \cite[Theorems 5.1 \& 6.1]{ByrNS94}.
Incidentally, since the {\small DFP} (Davidon-Fletcher-Powell) update is dual to the \bfgs by interchanging
$ \k{B} \leftrightarrow \k{H} $ and $ \k{s} \leftrightarrow \k{y} $ in eq. \eqref{eq:recIBFGS}, also
the {\small DFP} representation is known. Another well known rank-2 update for the direct Hessian approximation
is the Powell-Symmetric-Broyden (PSB) formula
\begin{equation}
    \label{eq:recPSB}
    \ko{B} = \k{B} + \frac{ ( \k{y} - \k{B} \k{s} )\kt{s} + \k{s}( \k{y} - \k{B} \k{s} )^T   }{\kt{s} \k{s} }
            - \frac{ ( \k{y} - \k{B} \k{s} )^T \k{s} }{ (\kt{s} \k{s})^2 } \k{s} \kt{s}
\end{equation}
For the {\small PSB} recursion, the compact representation has been discovered recently by Kanzow and Steck \cite{kanzowS23}.
Further compact representations for the Broyden class of updates have been developed in DeGuchy, Erway and Marcia \cite{DeGuchyEM18}. The compact
representation of structured \bfgs is proposed in Brust et \emph{al.} \cite{BrustDLP21} and for the multipoint-symmetric
secant matrix a representation has been developed in Brudakov et \emph{al.} \cite{BurMarPil02} and \cite{Brust18}. 
Because each of the recursive update formulas is unique (with its own advantages), the compact representations,
when they exist are also specific to the formula. More universally, in Dennis and Mor\'{e} \cite{DennisM77},
two general rank-2 update formulas are proposed that as special cases include the \bfgs formula (hence {\small DFP}),
the {\small PSB} update and the multipoint symmetric secant matrix among others. For \emph{arbitrary} vectors $ \k{v} \in \mathbb{R}^d $
and $ \k{c} \in \mathbb{R}^d $, so long $ \kt{v} \k{d} \ne 0  $ and $ \kt{c} \k{s} \ne 0 $ a general rank-2 formula
for the inverse is given by \cite[eq. 7.24]{DennisM77}
\begin{equation}
    \label{eq:recRk2H}
    \ko{H} = \k{H} + \frac{ ( \k{s} - \k{H} \k{y} )\kt{v} + \k{v}( \k{s} - \k{H} \k{y} )^T   }{\kt{v} \k{y} }
            - \frac{ ( \k{s} - \k{H} \k{y} )^T \k{y} }{ (\kt{v} \k{y})^2 } \k{v} \kt{v}
\end{equation}
and a general formula for the direct Hessian estimate is \cite[eq. 7.9]{DennisM77}
\begin{equation}
    \label{eq:recRk2B}
    \ko{B} = \k{B} + \frac{ ( \k{y} - \k{B} \k{s} )\kt{c} + \k{c}( \k{y} - \k{B} \k{s} )^T   }{\kt{c} \k{s} }
            - \frac{ ( \k{y} - \k{B} \k{s} )^T \k{s} }{ (\kt{c} \k{s})^2 } \k{c} \kt{c}
\end{equation}
It is straightforward to see that when $ \k{c} = \k{s} $ in eq. \eqref{eq:recRk2B} then this update corresponds
to the {\small PSB} formula in \eqref{eq:recPSB}. 
Moreover, it is known that when $ \k{v} = \k{s} $ in \eqref{eq:recRk2H} then the update
reduces to the inverse \bfgs formula from eq. \eqref{eq:recIBFGS}.
However, because the updates are parametrized by arbitrary vectors one can easily develop new methods
(like, for instance stochastic formulas), by replacing $\k{v}$ or $ \k{c} $ by other vectors. Nonetheless, to make
the formulas usable for large problems more effective representations than the dense formulas in 
\ref{eq:recRk2H} and \ref{eq:recRk2B} are needed. Similar to eqs. \eqref{eq:SYD} one can define the matrices
\begin{equation}
    \label{eq:VC}
     \k{V} = \bmat{ \subsc{v}{0} & \cdots & \subsc{v}{k-1} }, \quad \k{C} = \bmat{ \subsc{c}{0} & \cdots & \subsc{c}{k-1} }.  
\end{equation}

\subsection{Contributions}
\label{sec:contrib}
This work develops the compact representations of the dense matrix recurrences in equations \ref{eq:recRk2H} and \ref{eq:recRk2B}.
The representations enable effective limited-memory methods by storing only a small subset of previous vectors.
In particular, because the update formulas are defined in terms of arbitrary vectors (i.e, $ \k{v}$ and $\k{c}$), the
representations enable straightforward development of new methods, just by judiciously choosing particular 
sets of vectors in the representation. In this way, for instance, stochastic methods can be derived
by replacing deterministic quantities by random vectors. Further, we demonstrate how the representations
yield efficient eigenvalue computations, which make them viable for line-search as well as trust-region
optimization strategies.

\section{Compact Representations}
\label{sec:compactProof}
To develop the compact representation of eqs. \eqref{eq:recRk2H} and \eqref{eq:recRk2B} we use the subsequent notation:
For some matrix $ {X} \in \mathbb{R}^{d \times k}$ and $ {Y} \in \mathbb{R}^{ d \times k} $ we decompose the
product $ {X}^T {Y} $ as 
\begin{equation}
\label{eq:notatdecomp}
{X}^T {Y} = L^{XY} + R^{XY}, \quad \text{ and } \quad \text{diag}(\t{X}{Y}) = \ex{D}{XY}, 
\end{equation} 
where $ \ex{L}{XY} $ is the strictly lower triangular part and $ \ex{R}{XY} $ is the upper triangular part
(including the diagonal). These decompositions correspond to the element-wise definitions from eqs. \eqref{eq:SYD} and \eqref{eq:LR}, by
generalizing the $ \{ \subsc{s}{i} \} $'s with $ \{ \subsc{x}{i} \} $'s and using the $ \{ \subsc{y}{i} \} $'s. The only exception of the notation in \eqref{eq:notatdecomp} 
is when $ {X} = \k{S} $ and $ {Y} = \k{Y} $, in which case we omit the superscripts to be consistent with the notation in literature and simply write $ \kt{S} \k{Y} = \k{L} + \k{R} $,
with $\k{D} $ denoting the diagonal. But for any other value, e.g., $ {X} = \k{V}  $ and $ {Y} = \k{Y}  $ we write $ \kt{V} \k{Y} = \subex{L}{k}{VY} + \subex{R}{k}{VY} $. 
The result of the compact representation for the first rank-2 formula is stated as Theorem \ref{thm:compH}. 
\begin{theorem}
\label{thm:compH}
Applying the recursive update in eq. \eqref{eq:recRk2H} to a symmetric initialization
$ \subsc{H}{0} \in \mathbb{R}^{d \times d} $, with sequences $ \{ \subsc{s}{i} = \subsc{w}{i} - \subsc{w}{i-1}  \}_{i=0}^{k-1} $ and
$ \{ \subsc{y}{i} = \subsc{g}{i} - \subsc{g}{i-1}  \}_{i=0}^{k-1} $ and arbitrary vectors $ \{ v_i \}_{i=0}^{k-1} $ (so long $ \subsct{v}{i} \subsc{y}{i}  \ne 0 $) 
is equivalent to the compact representation
\begin{equation}
\label{eq:compthrm1}
    \k{H} = \subsc{H}{0} + 
            \bmat{ \k{V} \!& \k{S} \!-\! \subsc{H}{0} \k{Y}  }
            \bmat{0_{k \times k}    &   \subex{R}{k}{VY} \\ 
            \subext{R}{k}{VY}         & \k{R} \!+\! \kt{R} \!-\! ( \k{D} \!+\! \kt{Y} \subsc{H}{0} \k{Y}   )  }^{\!-\!1}
            \bmat{ \kt{V} \\ (\k{S} \!-\! \subsc{H}{0} \k{Y})^T   },
 \end{equation}
 where $\k{V}$, $ \k{S}, \k{Y}, \k{R} $ and $ \k{D} $ are defined in eqs. \eqref{eq:VC}, \eqref{eq:SYD}, \eqref{eq:LR} and $ \subex{R}{k}{VY}  $
 is the upper triangular part of $ \kt{V} \k{Y} $.
\end{theorem}
\begin{proof}
To simplify the expressions, we view the compact representation in eq. \eqref{eq:compthrm1} as
\begin{equation}
\label{eq:compactshort}
\k{H} = \subsc{H}{0} + \k{U} \k{M}^{-1} \kt{U},
\end{equation}
where $ \k{U} = \bmat{ \k{V} & \k{Z}  } $ and $ \k{Z}  = \k{S} - \subsc{H}{0} \k{Y} $ with columns  $ \subsc{z}{i} = \subsc{s}{i} - \subsc{H}{0} \subsc{y}{i}  $ and $ \k{M} $ being the middle matrix (e.g., the (1,2) block element is $ (\k{M})_{12} = \subex{R}{k}{VY} $ and $ (\k{M})_{11} = 0_{k \times k} $). The proof is by induction. 
We start with the base case, $k=1$. In this case $ \subex{R}{1}{VY} = \subsct{v}{0}\subsc{y}{0} $, $ \subsc{D}{1} = \subsct{s}{0} \subsc{y}{0}   $
and $ \subsc{R}{1} = \subsct{s}{0} \subsc{y}{0}  $, and the middle matrix $ \subsc{M}{1} $ from eq. \eqref{eq:compactshort} becomes
\begin{equation*}
\subsc{M}{1}^{-1} =
    \bmat{  0 & \subsct{v}{0}\subsc{y}{0} \\ 
    \subsct{y}{0}\subsc{v}{0} & \subsct{s}{0}\subsc{y}{0} - \subsct{y}{0} \subsc{H}{0} \subsc{y}{0}  }^{-1}
    =
    \frac{-1}{(\subsct{v}{0} \subsc{y}{0} )^2}
    \bmat{ \subsct{s}{0}\subsc{y}{0} - \subsct{y}{0} \subsc{H}{0} \subsc{y}{0} & -\subsct{v}{0}\subsc{y}{0} \\ 
    -\subsct{y}{0}\subsc{v}{0} & 0  }
\end{equation*}
Using the latter expression for $\subsc{M}{1}^{-1}$ and $ \subsc{U}{1} = \bmat{ \subsc{v}{0} & \subsc{z}{0}  }  $, the product $ \subsc{U}{1} \subsc{M}{1}^{-1} \subsct{U}{1}  $ is
\begin{equation*}
    \subsc{U}{1} \subsc{M}{1}^{-1} \subsct{U}{1}  = - \frac{\subsct{s}{0} \subsc{y}{0} - 
        \subsct{y}{0} \subsc{H}{0} \subsc{y}{0} }{ (\subsct{y}{0} \subsc{v}{0} )^2 } \subsc{v}{0}\subsct{v}{0} + 
        \frac{\subsc{v}{0} \subsct{z}{0} }{\subsct{y}{0} \subsc{v}{0}} +  \frac{\subsc{z}{0} \subsct{v}{0} }{\subsct{y}{0} \subsc{v}{0}}.
\end{equation*}
Substituting $ \subsc{z}{0} = \subsc{s}{0} - \subsc{H}{0} \subsc{y}{0}   $ and adding the initial matrix, one can see that 
\begin{equation*}
\subsc{H}{1} = \subsc{H}{0} + \subsc{U}{1} \subsc{M}{1}^{-1} \subsct{U}{1} = 
 \subsc{H}{0} - \frac{\subsct{s}{0} \subsc{y}{0} - 
        \subsct{y}{0} \subsc{H}{0} \subsc{y}{0} }{ (\subsct{y}{0} \subsc{v}{0} )^2 } \subsc{v}{0}\subsct{v}{0} + 
        \frac{\subsc{v}{0} (\subsc{s}{0} - \subsc{H}{0} \subsc{y}{0})^T }{\subsct{y}{0} \subsc{v}{0}} +  \frac{ ( \subsc{s}{0} - \subsc{H}{0} \subsc{y}{0} ) \subsct{v}{0} }{\subsct{y}{0} \subsc{v}{0}}
\end{equation*}
The last equality for $ \subsc{H}{1} $ is the same as the recursive update formula in eq. \eqref{eq:recRk2H}, and therefore the
compact formula holds in the base case. 
Assuming that the representation in eq. \eqref{eq:compactshort} is true for some $k \ge 1$, we now show that applying the rank-2 update in
eq. \eqref{eq:recRk2H} yields the representation at index $k+1$.
We start with the vector $ \k{s} - \k{H} \k{y}   $: 
\begin{equation}
    \label{eq:sHy}
    \k{s} - \k{H} \k{y} = \k{s} - ( \subsc{H}{0} + \k{U} \k{M}^{-1} \kt{U} ) \k{y} 
                        = \k{s} - \subsc{H}{0} \k{y} - \k{U} \k{M}^{-1} \kt{U}  \k{y} 
                        = \k{z} - \k{U} \k{M}^{-1} \kt{U}  \k{y},
\end{equation}
where $ \k{z} = \k{s} - \subsc{H}{0} \k{y}    $. Next define the scalars
\begin{equation}
    \label{eq:scalars}
    \k{\rho} = \kt{v} \k{y}, \quad \k{\beta} = \kt{y} \k{z}, \quad 
    \k{\theta} = \kt{y}( \k{s} - \k{H} \k{y}  ) = \k{\beta} - \kt{y} \k{U} \k{M}^{-1} \kt{U} \k{y}. 
\end{equation}
We now substitute eqs. \eqref{eq:sHy}, \eqref{eq:scalars} and \eqref{eq:compactshort} into the recursive formula eq. \eqref{eq:recRk2H} 
in order to rewrite it in terms of the compact representation
\begin{align}
\ko{H} \!&=\! \k{H} + \frac{ ( \k{s} - \k{H} \k{y} )\kt{v} + \k{v}( \k{s} - \k{H} \k{y} )^T   }{\kt{v} \k{y} }
            - \frac{ ( \k{s} - \k{H} \k{y} )^T \k{y} }{ (\kt{v} \k{y})^2 } \k{v} \kt{v} \nonumber  \\
        &= \k{H} + \frac{( \k{z} \!-\! \k{U} \k{M}^{-1} \kt{U}  \k{y} )\kt{v} + \k{v}(\k{z} \!-\! \k{U} \k{M}^{-1} \kt{U}  \k{y})^T }{\kt{v} \k{y}}  
        \!-\! \frac{ ( \k{s} \!-\! \k{H} \k{y} )^T \k{y} }{ (\kt{v} \k{y})^2 } \k{v} \kt{v} \nonumber \\
        &= \k{H} + \frac{1}{\k{\rho}}
        \bigg( 
            \bmat{ \k{U} \!& \k{v} \!& \k{z} \!}
            \left[
            \begin{array}{ c c  c }
                0_{2k \times 2k} & \!-\! \k{M}^{\!-\!1} \kt{U} \k{y} \!\! \! & 0_{2k \times 1}  \\
                \! \!-\!  \kt{y} \k{U} \k{M}^{\!-\!1} \!\! & 0 & 1 \\
                0_{1 \times 2k} \! \! & 1 & 0
            \end{array}
            \right]
            \bmat{ \kt{U} \\ \kt{v} \\ \kt{z} }
        \bigg)
        \!-\! \frac{\k{\theta}}{\k{\rho}^2} \k{v} \kt{v} \nonumber \\
        \label{eq:compblock}
        &= \subsc{H}{0} +
            \bmat{ \k{U} & \k{v} & \k{z} }
            \left[
            \begin{array}{ c c   c }
                \k{M}^{-1} & - \frac{1}{\k{\rho}}\k{M}^{-1} \kt{U} \k{y} & 0_{2k \times 1} \\
                -  \frac{1}{\k{\rho}}\kt{y} \k{U} \k{M}^{-1} & - \frac{\k{\theta}}{\k{\rho}^2} & \frac{1}{\k{\rho}} \\
                0_{1 \times 2k} & \frac{1}{\k{\rho}} & 0
            \end{array}
            \right]
            \bmat{ \kt{U} \\ \kt{v} \\ \kt{z} }
\end{align}
Next we compute an inverse representation of the block $ 3 \times 3 $ middle matrix in \eqref{eq:compblock} by inverting its 
upper $(2k+1) \times (2k+1) $ block first 
\begin{equation*}
\bmat{ \k{M}^{-1} & - \frac{1}{\k{\rho}}\k{M}^{-1} \kt{U} \k{y}  \\
                -  \frac{1}{\k{\rho}}\kt{y} \k{U} \k{M}^{-1} & - \frac{\k{\theta}}{\k{\rho}^2}} =
\bmat{ \k{M} - \frac{1}{\k{\beta}} \kt{U} \k{y} \kt{y} \k{U} & - \frac{\k{\rho}}{\k{\beta}}\kt{U} \k{y}  \\
                - \frac{\k{\rho}}{\k{\beta}} \kt{y} \k{U} & - \frac{\k{\rho}^2}{\k{\beta}}}^{-1}
\end{equation*}
Therefore, eq. \eqref{eq:compblock} becomes
\begin{align}
\ko{H} \!&=\!  \subsc{H}{0} +
            \bmat{ \k{U} & \k{v} & \k{z} }
            \left[
            \begin{array}{ c c |  c }
                \k{M}^{-1} & - \frac{1}{\k{\rho}}\k{M}^{-1} \kt{U} \k{y} & 0_{2k \times 1} \\
                -  \frac{1}{\k{\rho}}\kt{y} \k{U} \k{M}^{-1} & - \frac{\k{\theta}}{\k{\rho}^2} & \frac{1}{\k{\rho}} \\
                \hline 0_{1 \times 2k} & \frac{1}{\k{\rho}} & 0
            \end{array}
            \right]
            \bmat{ \kt{U} \\ \kt{v} \\ \kt{z} } \nonumber \\
        &= 
        \subsc{H}{0} +
            \bmat{ \k{U} & \k{v} & \k{z} }
            \left[
            \begin{array}{ c  |  c  }
                \bmat{ \k{M} \!-\! \frac{1}{\k{\beta}} \kt{U} \k{y} \kt{y} \k{U} & \!-\! \frac{\k{\rho}}{\k{\beta}}\kt{U} \k{y}  \\
                - \frac{\k{\rho}}{\k{\beta}} \kt{y} \k{U} & - \frac{\k{\rho}^2}{\k{\beta}}}^{\!-\!1} & 
                \begin{array}{c}
                    0_{2k \times 1} \\
                    \frac{1}{\k{\rho}}
                \end{array} \\
                \hline \begin{array}{c c}
                    0_{1 \times 2k} \phantom{- \frac{1}{\k{\beta}} \kt{U} \k{y} }&
                    \phantom{\kt{y}} \frac{1}{\k{\rho}}
                \end{array} & 0
            \end{array}
            \right]
            \bmat{ \kt{U} \\ \kt{v} \\ \kt{z} } \nonumber \\
        &= 
        \label{eq:compLarge}
        \subsc{H}{0} +
            \bmat{ \k{U} & \k{v} & \k{z} }
            \bmat{ \k{M} & 0_{2k \times 1} & \kt{U} \k{y} \\
                   0_{1 \times 2k} & 0 & \k{\rho} \\ 
                   \kt{y} \k{U} & \k{\rho} & \k{\beta} }^{-1}
            \bmat{ \kt{U} \\ \kt{v} \\ \kt{z} } 
\end{align}
Next we recall that $ \k{U} = \bmat{ \k{V} & \k{Z}  }  $, and decompose $ \k{M} $ into $ 2 \times 2 $ blocks,
with $ (\k{M})_{ij}, 1 \le i,j \le 2 $. 
Hence, the expression from
\eqref{eq:compLarge} is
\begin{align}
        \ko{H} 
        &=         
        \subsc{H}{0} +
            \bmat{ \k{U} & \k{v} & \k{z} }
            \bmat{ \k{M} & 0_{2k \times 1} & \kt{U} \k{y} \\
                   0_{1 \times 2k} & 0 & \k{\rho} \\ 
                   \kt{y} \k{U} & \k{\rho} & \k{\beta} }^{-1}
            \bmat{ \kt{U} \\ \kt{v} \\ \kt{z} } \nonumber \\
        &= \subsc{H}{0} +
            \bmat{ \k{V} & \k{Z} &  \k{v} & \k{z} }
            \bmat{ (\k{M})_{11} & (\k{M})_{12} &  0_{k \times 1} & \kt{V} \k{y} \\
                    (\k{M})_{21} & (\k{M})_{22} &  0_{k \times 1} & \kt{Z} \k{y} \\
                   0_{1 \times k} & 0_{1 \times k} &  0 & \k{\rho} \\ 
                   \kt{y} \k{V} & \kt{y} \k{Z} & \k{\rho} & \k{\beta} }^{-1}
                   \bmat{ \kt{V} \\ \kt{Z} \\ \kt{v} \\ \kt{z} } \nonumber \\
        &= \label{eq:compPatt} \subsc{H}{0} +
            \bmat{ \k{V} &  \k{v}  &  \k{Z} & \k{z} }
            \bmat{ (\k{M})_{11} & 0_{k \times 1} &   (\k{M})_{12} & \kt{V} \k{y} \\
                   0_{1 \times k}  & 0  &  0_{1 \times k}  & \k{\rho} \\
                    (\k{M})_{21} & 0_{k \times 1}  &  (\k{M})_{22} & \kt{Z} \k{y}  \\ 
                   \kt{y} \k{V} &  \k{\rho}& \kt{y} \k{Z} & \k{\beta} }^{-1}
                   \bmat{ \kt{V} \\ \kt{v} \\ \kt{Z} \\ \kt{z} }  
\end{align}
Since $ \k{Z} = \k{S} - \subsc{H}{0}\k{Y}   $, $\k{z} = \k{s} - \subsc{H}{0}\k{y} $,
$ \k{\rho} = \kt{v} \k{y}  $ and $ \k{\beta} = \kt{z} \k{y} = (\k{s} - \subsc{H}{0}\k{y})^T \k{y} $
we now verify the form of the matrices for $k+1$
\begin{align*}
    \bmat{ \k{V} & \k{v} } &= \ko{V} \\
    \bmat{ \k{Z} & \k{z} } &= \ko{Z} = \bmat{ \k{S} - \subsc{H}{0} \k{Y} & \k{s} - \subsc{H}{0} \k{y}  } = \ko{S} - \subsc{H}{0} \ko{Y}  \\
    \bmat{ (\k{M})_{11} & 0_{k \times 1} \\ 0_{1 \times k} & 0 } &= \bmat{ 0_{k \times k} & 0_{k \times 1} \\ 0_{1 \times k} & 0 } = 0_{(k+1) \times (k+1) } \\
    \bmat{ (\k{M})_{21} & 0_{k \times 1} \\ \kt{y} \k{V} & \k{\rho} } &= \bmat{ (\subex{R}{k}{VY})^T &  \\  \kt{y} \k{V} & \kt{y} \k{v} } 
    = (\subex{R}{k+1}{VY})^T \\
    \bmat{ (\k{M})_{12} & \kt{V} \k{y} \\ 0_{1\times k}  & \k{\rho} } &= \bmat{ (\subex{R}{k}{VY}) & \kt{V}\k{y}  \\   & \kt{y} \k{v} } 
    = \subex{R}{k+1}{VY} \\
    \bmat{ (\k{M})_{22} & \kt{Z} \k{y} \\ \kt{y} \k{Z}   & \k{\beta} } &= \bmat{\!\k{R}\!\!+\!\!\kt{R}\!\!-\!\!(\k{D}\!\!+\!\!\kt{Y} \subsc{H}{0} \k{Y} )   \!\!\!& \!(\k{S}\!\!-\!\!\subsc{H}{0}\k{Y} )^T\k{y}\!  \\  \kt{y}(\k{S}\!\!-\!\!\subsc{H}{0}\k{Y}  )  & \kt{y}( \k{s}\!\!-\!\!\subsc{H}{0} \k{y}   ) } \\
    & = R_{k\!+\!1}\!\!+\!\!R^T_{k\!+\!1}\!\!-\!\!(D_{k\!+\!1}\!\!+\!\!Y^T_{k\!+\!1} \subsc{H}{0} Y_{k\!+\!1} ) 
\end{align*}
Therefore, from eq. \eqref{eq:compPatt} and the matrices at $k+1$ the compact representation in eq. \eqref{eq:compactshort} becomes
\begin{equation*}
    \ko{H} = \subsc{H}{0} + \ko{U} \ko{M}^{-1} \kot{U},
\end{equation*}
with $ \ko{U} = \bmat{ \ko{V} & \ko{Z} } $ and the corresponding blocks for $ \ko{M} $. The explict components
are given in eq. \eqref{eq:compthrm1}, completing the induction.
\end{proof}
We note a few further observations about the compact representation in eq. \eqref{eq:compthrm1}. In a direct implementation
one stores and updates a few matrices that grow with $k \ge 1$, and uses a constant initialization. We will describe the updating techniques
in more detail in \cref{sec:lm}, but focus here on general memory properties. Suppose that
$ \subsc{H}{0} \in \mathbb{R}^{d \times d} $ is a constant diagonal, with $d$ nonzeros. Moreover, suppose one stores
and updates $ \k{V} \in \mathbb{R}^{d \times k} $, $ \k{Z} \in \mathbb{R}^{d \times k} $ (where $ z_i = s_i - H_0 y_i $),
and the triangular matrices $ \subex{R}{k}{VY} \in \mathbb{R}^{ k \times k } $ and $ \subex{R}{k}{ZY} \in \mathbb{R}^{ k \times k } $.
Storing these quantities accounts for all terms of the compact
representation, because $ \subex{D}{k}{ZT} = \text{diag}(\subex{R}{k}{ZY})  $ and 
$ \k{R} + \kt{R} - (\k{D} + \kt{Y} \subsc{H}{0} \k{Y}) = \subex{R}{k}{ZY} + \subext{R}{k}{ZY} - \subex{D}{k}{ZY}  $.
The memory of eq. \eqref{eq:compthrm1} with this storage scheme is
\begin{equation}
    \label{eq:memk}
    d + 2dk + 2 \left( \frac{k(k+1)}{2} \right) = \mathcal{O}(2(kd + k^2/2))
\end{equation}
For large and difficult problems, where $d$ is large and many iterations $k$ are computed
the memory complexity of eq. \eqref{eq:memk} is not practical. However, for such situations,
a limited-memory technique can be efficiently implemented with the compact representation.
The initialization is typically chosen as a multiple of the identity,
which is updated each iteration $ \subsup{H}{k}{0} = \k{\gamma} I  $, $ \k{\gamma} \in \mathbb{R}  $.
Because the initialization changes every iteration $ \k{V} $, $ \k{S} $ and $ \k{Y} $ have to be stored
separately. However, for a small constant memory parameter $ l \ll d $ (say, $l=5$), the matrices
$ \k{V} $, $ \k{S} $ and $ \k{Y} $ are stored and are defined only by the $ l $ most recent updates, 
hence each being of size $ ld $. 
The limited memory sizes of $ \subex{R}{k}{VY} $, $ \k{R} $ and $ \subex{Y}{k}{YY} $ 
(the upper triangular part of $ \kt{Y} \k{Y} $ ) are $ \frac{(l+1)l}{2} $, respectively.
Therefore, the storage of the limited memory compact representation in eq. \eqref{eq:compthrm1} is
independent of $k$ and is given by
\begin{equation}
    \label{eq:meml}
1 + 3ld + 3 \left( \frac{l(l+1)}{2} \right) = \mathcal{O}(3(ld + l^2/2))
\end{equation}
Since $ l \ll d $ for most practical applications, the memory requirement is linear 
(a small constant multiple) in the size of the problem, i.e.,  $ \mathcal{O}(3ld) $.
This setting corresponds to the factored form in the right hand side of Fig. \ref{fig:compact}. 
In addition, note that the middle matrix in the compact representation eq. \eqref{eq:compthrm1}
can also be expressed with its explicit inverse
\begin{equation}
    \label{eq:invM}
    \k{M}^{-1} =
    \bmat{ - (\subex{R}{k}{VY})^{-T} \big( \k{R} + \kt{R} - (\k{D} + \kt{Y} \subsc{H}{0} \k{Y}  ) \big) (\subex{R}{k}{VY})^{-1}   & (\subex{R}{k}{VY})^{-T} \\ 
            (\subex{R}{k}{VY})^{-1} &  0_{k \times k} },
\end{equation}
\begin{equation}
    \label{eq:M}
    \k{M} =
    \bmat{ 0_{k \times k} & \subex{R}{k}{VY} \\ 
            (\subex{R}{k}{VY})^T & \k{R} + \kt{R} - (\k{D} + \kt{Y} \subsc{H}{0} \k{Y}  )  }
\end{equation}
(To verify these identities, compute e.g., $ \k{M}^{-1} \k{M} = I $). The explicit
inverse in eq. \eqref{eq:invM} shows that solves with $ \k{M}^{-1} $ can be computed efficiently
with solves of two triangular matrices $ (\subex{R}{k}{VY})^T $ and $ \subex{R}{k}{VY} $,
and overall multiplication complexity of $ \mathcal{O}( l^2 )   $. Furthermore,
because $ (\subex{R}{k}{VY} )_{ii} = \subsct{v}{i} \subsc{y}{i} $ and $ \subex{R}{k}{VY}  $
is upper triangular, eq. \eqref{eq:invM} implies that $ \subsct{v}{i} \subsc{y}{i} \ne 0 $
is a necessary and sufficient condition for the existence of the compact representation in
eq. \eqref{eq:compthrm1}.
As a useful by-product we can use the derivations for Theorem \ref{thm:compH} to deduce the compact representation
of the recursive update formula for the direct Hessian approximation $ \k{B} $ in eq. \eqref{eq:recRk2B}.
The following theorem is symmetric to the previous result:
\begin{theorem}
\label{thm:compB}
Applying the recursive update in eq. \eqref{eq:recRk2B} to a symmetric initialization
$ \subsc{B}{0} \in \mathbb{R}^{d \times d} $, with sequences $ \{ \subsc{y}{i} = \subsc{g}{i} - \subsc{g}{i-1}  \}_{i=0}^{k-1} $ and
$ \{ \subsc{s}{i} = \subsc{w}{i} - \subsc{w}{i-1}  \}_{i=0}^{k-1} $ and arbitrary vectors $ \{ c_i \}_{i=0}^{k-1} $ (so long $ \subsct{c}{i} \subsc{s}{i}  \ne 0 $) 
is equivalent to the compact representation
\begin{equation}
\label{eq:compthrm2}
    \k{B} = \subsc{B}{0} + 
            \bmat{ \k{C} & \k{Y} \!-\! \subsc{B}{0} \k{S}  }
            \bmat{0_{k \times k}    &   \subex{R}{k}{CS} \\ 
            \subext{R}{k}{CS}         & \k{R} \!+\! \kt{R} \!-\! ( \k{D} \!+\! \kt{S} \subsc{B}{0} \k{S}   )  }^{\!-\!1}
            \bmat{ \kt{C} \\ (\k{Y} \!-\! \subsc{B}{0} \k{S})^T   },
 \end{equation}
 where $ \k{C},  \k{S}, \k{Y}, \k{R} $ and $ \k{D} $ are defined in eqs. \eqref{eq:VC}, \eqref{eq:SYD}, \eqref{eq:LR} and $ \subex{R}{k}{CS}  $
 is the upper triangular part of $ \kt{C} \k{S} $.
\end{theorem}
\begin{proof}
Observe that the recursive update for $ \ko{B} $ in eq. \eqref{eq:recRk2B} can be obtained from eq. \eqref{eq:recRk2H} by interchanging
$ \k{H} \leftrightarrow \k{B} $, $ \k{y} \leftrightarrow \k{s} $ and $ \k{v} \leftrightarrow \k{c} $. We now apply the same changes
to the compact representation in eq. \eqref{eq:compthrm1}. Specifically, $ \subsc{H}{0} \leftrightarrow \subsc{B}{0} $, 
$ \k{Y} \leftrightarrow \k{S} $ and $ \k{V} \leftrightarrow \k{C} $ gives the representation in eq. \eqref{eq:compthrm2}.
\end{proof}
We note that when $ \k{C} = \k{S} $ in \cref{eq:compthrm2} then this representation reduces to the PSB compact representation.

\subsection{Implications}
\label{sec:impl}
We develop further consequences of the compact representations in this section. Initially we focus
on the inverse representation from eq. \eqref{eq:compthrm1}, since many results carry over to the
direct factorization in eq. \eqref{eq:compthrm2} by symmetrically interchanging variables.
Moreover, most proofs are given in the appendix in order to avoid distraction from the main observations.
First, for limited-memory implementations it is standard to use a multiple of the identity initialization
that dynamically varies for every iteration, $ \subsup{H}{k}{0} = \k{\gamma} I  $. In this situation,
and whenever the initialization dynamically changes, 
one has to store $ \k{S} $ and $ \k{Y} $ separately in order to define $ \k{S} - \subsup{H}{k}{0} \k{Y} $.
However, it can be desirable to not form $ \k{S} - \subsup{H}{k}{0} \k{Y} $ explicitly (see e.g, \cite{kanzowS23}, \cite[Theorem 2]{brustEM24}).
The main approach to achieve this, is by expressing $ \k{S} - \subsup{H}{k}{0} \k{Y} $ as a product
\begin{equation*}
    \k{S} - \subsup{H}{k}{0} \k{Y} = \bmat{ \k{S} &  \subsup{H}{k}{0} \k{Y} } \bmat{ I \\ - I }
\end{equation*}
In a Corollary to Theorem \ref{eq:compthrm1} we describe a formulation that decouples $ \k{S} - \subsup{H}{k}{0} \k{Y}  $,
by storing $ \k{S} $ and $ \subsup{H}{k}{0} \k{Y} $ with a non-constant initialization $ \subsup{H}{k}{0} $.
\begin{corollary}
\label{cor:alt}
An alternative to the compact representation in Theorem \ref{eq:compthrm1} with $ \subsc{H}{0} = \subsup{H}{k}{0} $ and
$ \k{M}^{-1}  $ in \eqref{eq:invM} is given by
\begin{equation}
    \label{eq:compHY}
        \k{H} = \subsup{H}{k}{0} + 
        \bmat{\k{V} & \k{S}~| & \subsup{H}{k}{0} \k{Y} }
        \left[
        \begin{array}{c c | c}
            & &  \\
            &  \bigg[ \k{M}^{-1} \bigg]_{\phantom{A}}  & \bmat { - (\subex{R}{k}{VY})^{\!-\!T\!} \\ 0 }  \\
            \hline & \bmat { - (\subex{R}{k}{VY})^{\!-\!1\!} & 0 }   & 0
        \end{array}
        \right]
        \left[
            \begin{array}{c}
                \kt{V} \\
                \kt{S} \\
                \kt{Y} \subsup{H}{k}{0}
            \end{array}
        \right]
\end{equation}
\end{corollary}
\begin{proof}
The proof is in appendix \ref{app:A}.
\end{proof}
Corollary \ref{cor:alt} separates $ \k{S} $ from $ \subsc{H}{0} \k{Y} $, but it also establishes a connection
to the inverse \bfgs compact representation from eq. \eqref{eq:compactIBFGS}, because it also stores 
$ \k{S} $ and $ \subsc{H}{0} \k{Y} $ as the \bfgs representation does. In fact, when $ \k{V} = \k{S} $
the representation in Corollary \ref{cor:alt} (and hence Theorem \ref{eq:compthrm1}) is equivalent to the 
compact representation of the \bfgs formula:
\begin{corollary}
\label{cor:reductbfgs}
By choosing $ \k{V} = \k{S} $ in eqs. \eqref{eq:compHY} or equivalently \eqref{eq:compthrm1} the representation reduces
to the compact \bfgs formula in eq. \eqref{eq:compactIBFGS}.
\end{corollary}
\begin{proof}
The proof is in the supplemental materials . 
\end{proof}

Recall, since $\k{V}$ can be determined by choice other representations can therefore be designed by
substituting for this matrix. 
For instance, when $ \k{v} = \k{y} $ in the recursive update eq. \eqref{eq:recRk2H} then this formula
is known as Greenstadt's update. However, to the best of our knowledge, no compact representation
for this update has been discovered yet. Simply replacing $V_k=Y_k$ in \eqref{eq:compactIBFGS} gives this representation.  In corollary \ref{cor:compGreen} we describe another
compact representation for this recursion for the case that $H_0$ is a multiple
of the identity $ H_0 = \gamma_k I  $.

\begin{corollary}
\label{cor:compGreen}
The compact representation of the recursive update \eqref{eq:recRk2H} with $\k{v} = \k{y} $,
and $ H_0 = \gamma_k I $ also
known as Greenstadt's formula \cite[Section 7.3]{DennisM77}, is given by 
\begin{equation}
    \label{eq:compGreenstadt}
    \k{H} = \subsc{H}{0} + \bmat{ \k{S} & \subsc{H}{0} \k{Y}   } \k{N}^{-1} \bmat{ \kt{S} \\  \kt{Y} \subsc{H}{0}   },  
\end{equation}
\begin{equation}
    \label{eq:mcompGreenstadt}
    \k{N} = 
    \bmat{  \k{R} \!\!+\!\! \kt{R} \!\!-\!\! (\k{D} \!\!+\!\! \kt{Y} \subsc{H}{0} \k{Y}) 
        \!\!+\!\! \gamma_k\subex{R}{k}{YY} \!\!+\!\! \gamma_k\subext{R}{k}{YY}     &   \gamma_k\subext{R}{k}{YY} \\ 
            \gamma_k\subex{R}{k}{YY}         & 0_{k \times k}  },
\end{equation}
 where $ \k{S}, \k{Y}, \k{R} $ and $ \k{D} $ are defined in eqs. \eqref{eq:SYD} and \eqref{eq:LR} and $ \subex{R}{k}{YY}  $
 is the upper triangular part of $ \kt{Y} \k{Y} $.
\end{corollary}
\begin{proof}
The proof is in appendix \ref{app:B}.
\end{proof}

Because the compact formula in eq. \eqref{eq:mcompGreenstadt} is new we verify the validity of it by comparing
it to the recursive update from eq. \eqref{eq:recRk2H} in \cref{tab:greenstadtErrs}.

\begin{table}
\centering
\begin{tabular}{| c c c |}
\hline $k$ & Error 1 & Error 2 \\
\hline
1 	 &\texttt{5.63e-16} 	 &\texttt{5.17e-16} \\ 
2 	 &\texttt{1.17e-15} 	 &\texttt{9.15e-16} \\ 
3 	 &\texttt{6.88e-16} 	 &\texttt{1.12e-15} \\ 
4 	 &\texttt{1.2e-15} 	 &\texttt{1.28e-15} \\ 
5 	 &\texttt{1.75e-15} 	 &\texttt{1.47e-15} \\ 
6 	 &\texttt{1.64e-15} 	 &\texttt{1.69e-15} \\ 
7 	 &\texttt{2.74e-15} 	 &\texttt{2.06e-15} \\ 
8 	 &\texttt{3.6e-15} 	 &\texttt{2.58e-15} \\ 
\hline
\end{tabular}
\caption{Differences between the recursive rank-2 update eq. \eqref{eq:recRk2H}, $\k{H}^{\text{R}}$
with $ \k{v} = \k{y}  $ and the compact representation in \cref{cor:compGreen}, $\k{H}^{\text{C}}$.
Error 1 denotes the residual $\| \k{H}^{\text{C}} \subsc{y}{k-1} - \subsc{s}{k-1}   \|_2$, and 
Error 2 is the difference  $\| \k{H}^{\text{C}} - \k{H}^{\text{R}}   \|_F$}.
\label{tab:greenstadtErrs}
\end{table}

\subsection{Eigendecomposition}
\label{sec:eigendecomp}
For limited-memory methods it is common to use a multiple of the identity initialization $ \subsup{B}{k}{0} = (\subsup{H}{k}{0})^{-1} =  \k{\gamma}^{-1} I $
so that the compact representation (e.g., for the direct Hessian) can be viewed as 
\begin{equation}
    \label{eq:compeig}
    \k{B} = \subsup{B}{k}{0} + \k{J} \k{K}^{-1} \k{J}^T = \frac{1}{\k{\gamma}} I + \k{J} \k{K}^{-1} \k{J}^T,
\end{equation}
where $ \k{J} $ and $ \k{K} $ are specified by appropriate formulae (for instance,  eq. \eqref{eq:compthrm2}). Nonetheless, $ \k{J} $ is typically
very tall and skinny with dimension, say $ d \times 2l $  and $ l \ll d $. Hence, $ \k{K}  $ is a small symmetric square
of size $ 2l \times 2l $. We suppress the iteration index for the moment, and assume that $ J $ is of size $ d \times 2l $.
It is possible to exploit the representation in eq. \eqref{eq:compeig} in order to compute and \emph{implicit} eigendecomposition
with complexity that is linear in $d$. Suppose the ``thin'' QR factorization of $ J $ is $ J = QR $ (at about $ \mathcal{O}(4dl^2)  $ multiplications).
Then compute a small eigendecomposition of $ R K^{-1} R^T $  at $ \mathcal{O}(16l^3)  $ multiplications
\begin{equation*}
    R K^{-1} R^T = \widehat{P} \widehat{\Lambda} \widehat{P}^T,
\end{equation*}
where $ \widehat{P} \in \mathbb{R}^{2l \times 2l} $ is orthogonal and $ \widehat{\Lambda} \in \mathbb{R}^{2l \times 2l} $ is diagonal.
Define the thin orthonormal matrix $ \subsc{P}{(1)} = Q \widehat{P} \in \mathbb{R}^{d \times 2l}  $ and also its orthogonal
complement $ \subsc{P}{(2)} \in \mathbb{R}^{ d \times (d-2l) } $ (so that $ \subsc{P}{(2)}^T \subsc{P}{(1)} =0  $ and $ \subsc{P}{(2)}^T \subsc{P}{(2)} = I  $).
Note that the factors of  $ \subsc{P}{(1)} $, i.e., $ Q $ and $ \widehat{P} $ are explicitly computed, however the potentially
very large $ \subsc{P}{(2)} $ is never formed (it is only defined implicitly). Denote the eigenvalues corresponding to the 
eigenvectors in $ \subsc{P}{(1)} $ by
\begin{equation*}
    \lambda_i = \widehat{\lambda}_i + \frac{1}{\gamma}, \quad 1 \le i \le 2l,
\end{equation*}
and the remaining eigenvalues corresponding to the eigenvectors in $ \subsc{P}{(2)} $ by
\begin{equation*}
    \lambda_i = \frac{1}{\gamma}, \quad 2l+1 \le i \le d.
\end{equation*}
Representing the orthogonal matrix $ P = \bmat{ \subsc{P}{(1)} & \subsc{P}{(2)} } \in \mathbb{R}^{d \times d} $ and the diagonal matrix of $ \lambda_i $'s 
as $ \Lambda = \text{diag}(\lambda_1, \lambda_2,\ldots, \lambda_{d}) = \text{blkdiag}(\subsc{\Lambda}{(1)},\frac{1}{\gamma} I) \in \mathbb{R}^{d \times d } $,
the eigendecomposition of $B$ is
\begin{equation}
    \label{eq:eigen}
    B = \frac{1}{\gamma} I + J K J^T = \frac{1}{\gamma} I + QR K R^T Q^T = \frac{1}{\gamma} I +  Q \widehat{P} \widehat{\Lambda} \widehat{P}^T Q^T = P \Lambda P^T
\end{equation}
The factorization in eq. \eqref{eq:eigen} is implicit, because $ \subsc{P}{(2)} $ is never fully computed. Since,
$ \frac{1}{\gamma} $ is a repeated eigenvalue corresponding to the eigenspace of $ \subsc{P}{(2)} $ one can compute the scaled
projections  $ \frac{1}{\gamma} \subsc{P}{(2)} \subsc{P}{(2)}^T $ using the available $ \subsc{P}{(1)} = Q \widehat{P} $ only
\begin{equation*}
    \frac{1}{\gamma} \subsc{P}{(2)} \subsc{P}{(2)}^T = \frac{1}{\gamma} (I - \subsc{P}{(1)} \subsc{P}{(1)}^T ) = \frac{1}{\gamma} (I - Q Q^T).
\end{equation*}
This latter identity uses the orthogonality of $P$, that is $ I = P P^T = \subsc{P}{(1)} \subsc{P}{(1)}^T + \subsc{P}{(2)} \subsc{P}{(2)}^T   $.
For limited memory methods with large $d$ and small $l$, the main computational cost in computing the implicit eigendecomposition
of the compact representation is a thin QR factorization with linear complexity in the dimension of the problem. Therefore,
the eigendecomposition can be computed efficiently. Especially, for trust-region optimization methods computing the 
eigendecomposition is useful, because it enables effective shifting strategies that ensure positive definiteness of the matrix.

\subsection{Limited-Memory Updating}\label{sec:lm}
For large problems, limited-memory approaches store only a small number of 
vectors to define the representations. Depending on the initialization strategy, specifically
whether $ \subsup{H}{k}{0}$ varies between iterations or is constant
the matrices can be stored and updated in different ways (see the discussion following \cref{thm:compH}). 
We will describe some general techniques that apply to any initialization strategy in this section.
By setting the parameter $l \ll d$ limited-memory techniques enable inexpensive computations, and replace or insert one column at each
iteration in $ \k{Y} $, $ \k{S} $ and $ \k{V}$.
Let an underline below a matrix represent the 
matrix with its first column removed. That is, $ \underline{{S}}_k $ represents $ \k{S} $ without its first column.
With this notation, a column update of a matrix, say $ \k{S} $, by a vector $ \k{s} $ is defined as follows.
\begin{equation*}
	\text{colUpdate}\left(\k{S},\k{s} \right) \equiv
	\begin{cases}
		[\: \k{S} \: \k{s}\:  ] 						& \text{ if } k < l \\
		[\: \underline{S}_k \: \k{s}\:  ] 			& \text{ if } k \ge l. \\
	\end{cases}
\end{equation*}
This column update can be implemented efficiently, without copying large amounts of memory,
by appropriately updating the relevant index information. 
Certain matrix products can also be efficiently updated. As such, products like $ \k{S}^T \k{Y} $
do not have to be re-computed from scratch. In order to describe the matrix product updating mechanism,
let an overline above a matrix represent
the matrix with its first row removed. That is, $ \overline{\kt{S} \k{Y}} $ represents $ \kt{S} \k{Y} $ without its first row.
With this notation, a product update of, say $ \k{S}^T\k{Y} $, by matrices $ \k{S} $, $\k{Y} $
and vectors $ \k{s} $, $\k{y}$ is defined as:
\begin{equation*}
	\text{prodUpdate} \left( \k{S}^T\k{Y}, \k{S}, \k{Y}, \k{s}, \k{y} \right) \equiv
	\begin{cases}
		\left[
			\begin{array}{ c c }
				\k{S}^T\k{Y} 		& \k{S}^T\k{y} \\
				\k{s}^T\k{Y}	& \k{s}^T \k{y} 
			\end{array}
		\right] & \text{ if } k < l \vspace{0.1cm} \\		
		\left[
			\begin{array}{ c c }
				\left(\underline{\overline{\k{S}^T\k{Y}}}\right) 			& 	\underline{S}_k^T\k{y} \\
				\k{s}^T \underline{Y}_k					&	 \k{s}^T \k{y} 
			\end{array}
		\right] & \text{ if } k \ge l.\\
	\end{cases}
\end{equation*}
This product update can be implemented without recomputing potentially large multiplications, by storing
previous products and information about the column order in $ \k{S} $ and $ \k{Y} $. In particular, updating the matrix product can be done
by storing $ \k{S}^T \k{Y} $, $ \k{S}, \k{Y} $ and an appropriate vector of indices. 
Note that such a product update is computationally much more efficient, than recomputing the
product from scratch. Specifically, when $ l \le k $, the direct product $ \k{S}^T \k{Y} $ is
done at $ \mathcal{O}(l^2d) $ multiplications. However, an implementation of ``prodUpdate''  
does this update with $ \mathcal{O}(2ld) $ multiplications,
by reusing previous values represented by $\underline{\overline{\k{S}^T\k{Y}}}$. 
Moreover, when the product is symmetric, e.g. 
prodUpdate is invoked by e.g., $ \textnormal{prodUpdate}( \k{S}^T\k{S}, \k{S}, \k{S}, \k{s}, \k{s}) $, 
then $ \underline{S}^T_k \k{s} $ can be stored and reused in two places
(thus only one matrix-vector product is needed, instead of two). Updates to diagonal, lower or upper triangular
matrices can be described in a similar way. For instance, an update to the upper triangular matrix 
$ \k{R} $ can be computed via: $ \textnormal{prodUpdate}( \k{R}, \k{S}, 0, 0, \k{y}) $. These updating
techniques ensure that limited-memory computations retain their linear complexity with regards
to problem dimension.

\section{Numerical Experiments}
\label{sec:experiments}
In a set of numerical experiments we demonstrate the efficacy of the compact representations for a series
of data fitting tasks. A dedicated implementation of the corresponding algorithms is content for 
future research. The methods are implemented in {\small MATLAB} and {\small Python 3}
on a linux machine with intel 13th Gen Intel Core i9-13900KS (24 cores) processor and 128 GiB RAM.
All software is available in the public domain \url{https://github.com/johannesbrust/CR} 

\subsection{Eigenfactorization}
\label{sec:numexeigen}
\begin{figure}
\includegraphics[trim={0.1cm 0cm 0.0cm 0cm},clip,scale=0.68]{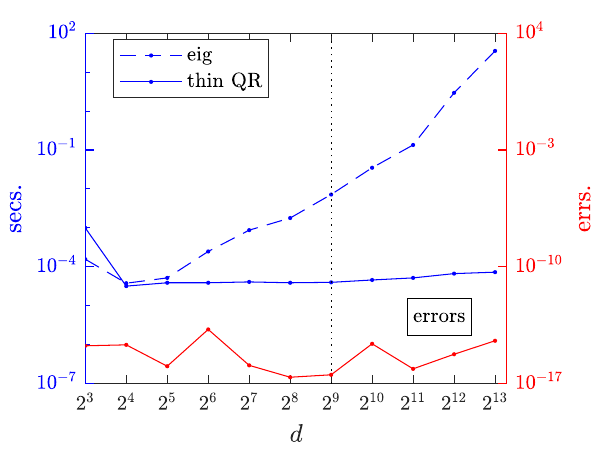}%
\includegraphics[trim={0.6cm 0cm 0.6cm 0cm},clip,scale=0.68]{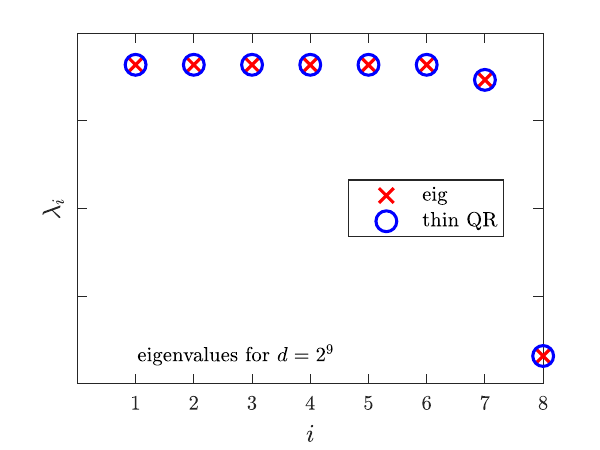}
\caption{Computing the eigenvalues of a compact representation in an optimization algorithm for the 
Rosenbrock function with $d=2^3,2^4\ldots,2^{13}$. Using \texttt{eig} \cite{MATLAB23a}
scales cubically, while a thin QR algorithm grows linearly with problem size (left figure blue axis). The magnitude of the errors remains low: 
$ \textnormal{error} = (\sum_{i=1}^d (\lambda^{\textnormal{eig}}_i -  \lambda^{\textnormal{qr}}_i)^2)^{\frac{1}{2}} / d   $
(left figure red axis). For $d=2^9$ the first 8 eigenvalues are
computed using \texttt{eig} and the proposed QR approach in the right hand figure.}
\label{fig:eig}
\end{figure}

This experiment demonstrates the scalability of 
an eidendecomposition with the proposed representations.
To generate relevant matrices, we apply an minimization algorithm combined with
the compact representation from eq. \eqref{eq:compactIBFGS} with $\k{V} = \k{S}$ to
the even Rosenbrock function
\begin{equation*}
    f(w) = \sum_{i=1}^{d/2} 100 (w^2_{2i-1} - w_{2i} )^2 + (w_{2i-1} - 1 )^2
\end{equation*}
We compute the eigenvalues of the compact representation using a thin QR factorization as described in
Sec. \ref{sec:eigendecomp}. The dimensions are $ d \in \{ 2^3, 2^4, \ldots, 2^{13} \} $, and the limited-memory parameter is $l=5$. At 
iteration $k=10$ the eigenfactorization is computed by \texttt{eig} \cite{MATLAB23a} and by the thin QR factorization.
We record the time for each of the approaches and the corresponding errors in \cref{fig:eig}. Since
trust-region algorithms may use eigenvalues to compute subproblem solutions, and the factorizations
scale favorably with problem dimensions, the compact representation appears well suited for large
trust-region strategies.

\subsection{Tensor fitting}
\label{sec:numextensor}

\begin{figure}
\includegraphics[trim={0.1cm 0cm 0.1cm 0cm},clip,scale=0.68]{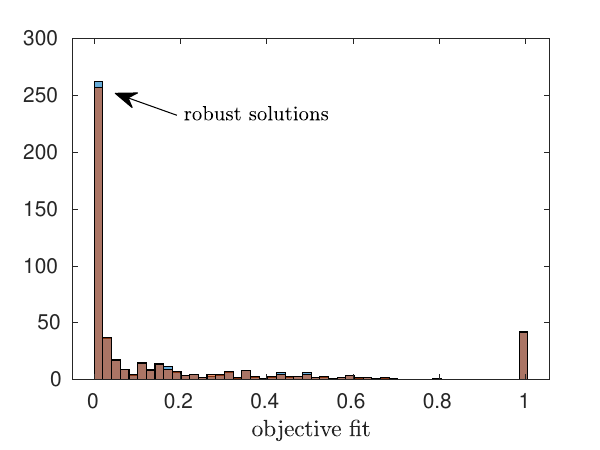}%
\includegraphics[trim={0.5cm 0cm 0.6cm 0cm},clip,scale=0.68]{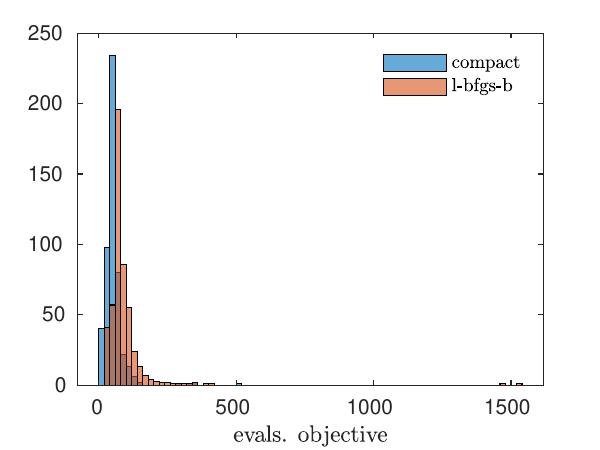}
\caption{The compact representation and algorithm \texttt{l-bfgs-b} are used to fit CP tensors.}
\label{fig:dist}
\end{figure}

In this experiment we use the compact representation to compute tensor factorizations. In particular, we compute the
Canonical polyadic decomposition (CP decomposition) with a given target rank $r$. The decomposition for 
a tensor in $ \mathbb{R}^{ d_1 \times \cdots \times d_m} $ is
\begin{equation}
\label{eq:tensor}
\mathcal{A} = \sum_{i=1}^r a_{d_1} \otimes a_{d_2} \otimes \cdots \otimes a_{d_m}
\end{equation}
This generalizes a low rank matrix approximation to higher order tensors. In order to fit the 
factorization to a given data tensor $ \widehat{\mathcal{A}} $ a nonlinear least-squares objective
is effective $ \min \| \widehat{\mathcal{A}} - \mathcal{A} \|_{F}  $. Because of the nonlinear form
of the factorization, however, the problem is typically nonconvex and multiple local solutions exist.
A solver of choice for
this fitting problem is \lbfgs \cite{AcarDK11}. We use the compact representation with $ \k{V} = \k{S}  $ 
in a strong Wolfe line-search to compute the tensor factorization. The limited-memory parameter is $ l=5 $
and the stopping condition is $ \| \k{g} \|_{\infty}  \le 1 \times 10^{-5}  $ for all solvers. We use the
{\small Tensor Toolbox (Sandia Natl Labs \& MathSci.ai)} \cite{TensorTBX} to generate the problems; the data tensors are of size
$ 250 \times 250 \times 250 $ and the target rank is $ r=2 $. The default method in the 
toolbox is {\small L-BFGS-B} \cite{ZhuByrdNocedal97} with a wrapper of the C implementation 
from \cite{BeckerLbfgs}. Five hundred tensor factorizations are solved for 
which we record the results. \cref{fig:dist} shows the distributions of the 
final fitted objective values and the number of function evaluations. Using the compact
solver results in a slightly higher frequency of lowest objective values (and therefore more robust
tensor reconstructions). Our compact implementation is effective in terms of total function evaluations
(\cref{fig:dist} right hand plot), which can be the main computational cost for large tensors. 

\subsection{A Multiclass model}
\label{sec:ml}

\begin{figure}
\includegraphics[trim={0.1cm 0cm 0.5cm 0cm},clip,scale=0.68]{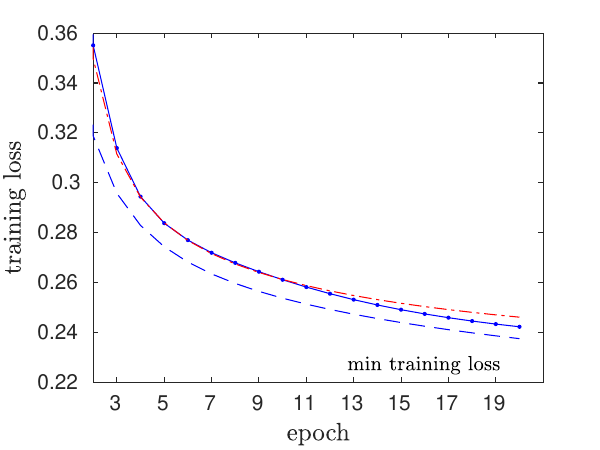}%
\includegraphics[trim={0.1cm 0cm 0.6cm 0cm},clip,scale=0.68]{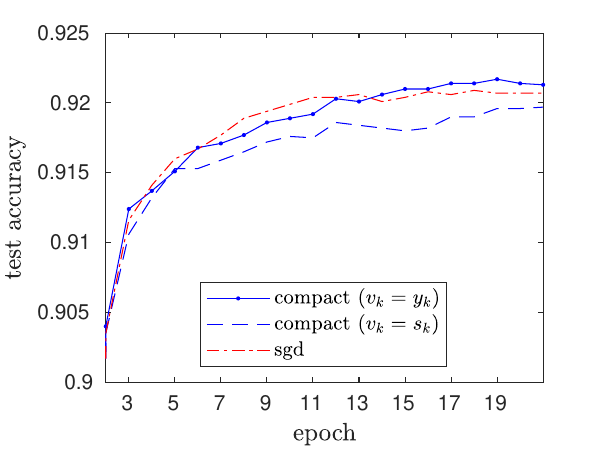}
\caption{Two compact representation algorithms and sgd \cite{robbinsM51} are used on a stochastic machine learning model.}
\label{fig:ml}
\end{figure}

In this experiment we implement a model, as well as stochastic minimization algorithms for a multiclass logistic regression of
the images in the MNIST dataset \cite{deng2012mnist}. In particular, the model should predict the correct label for a handwritten digit
as accurately as possible.  
Each of the images corresponds to an $ 28 \times 28 $ pixel array, $X_i$. Each label is a digit, $y_i \in \{0,1,2,\ldots,9 \}$.
The MNIST dataset contains $N=60,000$ images and labels as part of its training set. There are another $10,000$ pairs for a testset.
A conventional model first transforms the images into vectors $X_i \rightarrow x_i $ where each vector is in $\mathbb{R}^{28^2 \times 1} $
(e.g., by stacking the array columns onto each other) and then applies a parameter matrix $ W \in \mathbb{R}^{10 \times 28^2}$ to the image
$ W x_i $. In order for the model to use something that resembles probabilities one uses element-wise exponentiation to obtain nonnegative values. 
Specifically, the probability that the $i^{\text{th}}$ image has label $y_i$ is modeled by
\begin{equation*}
    \frac{ e^{[Wx_i]_{y_i +1}} }{ \sum_{j=1}^{10} e^{[Wx_i]_{j}} }
\end{equation*}
By applying the natural log to these probabilities and summing over all data points the loss and corresponding fitting problem is given by
\begin{equation}
    \label{eq:mainmin}
    \underset{ W }{ \textnormal{ minimize } } 
    \sum_{i=1}^N \ln \big( \sum_{j=1}^{10} e^{[Wx_i]_{j}}  \big) - [Wx_i]_{y_i +1}  
\end{equation}
The minimization problem in eq. \eqref{eq:mainmin} is of the form $ \text{min}  \sum_i f_i(W) $ with appropriate $f_i$. 
    In order to implement an (stochastic) algorithm one typically computes 
    a gradient with respect to the parameters represented as a vector. Thus we reparametrize the weight matrix
    as a vector $W \rightarrow w $ and then find the corresponding gradient $\nabla_w f_i$.
    Note that these gradients can also be used if a mini-batch method is used with subsampled gradients
    \begin{equation*}
        \frac{1}{N_i}\sum_{j=1}^{N_i} \nabla_w f_{i_j},  
    \end{equation*}
    and where $N_i$ is the size of the $i^{\text{th}}$ batch and the indices occur exactly once in the whole data set $ i_j \in \{1,2,\ldots N \} $.
    Using minibatch techniques results in stochastic problems, since every function evaluation is based on a (random) subset of the whole dataset.
    This means that traditional methods, such as line-search algorithms, are typically not advisable, because they rely on deterministic changes
    in function values. For a stochastic variation of the compact representation, we fix a constant learning rate (step size), $\alpha = \frac{1}{2}$
    at every iteration and update the iterates as $ \ko{w} = \k{w} + \alpha \k{p}  $. The step $ \k{p} $ is computed via the compact representation
    $ \k{p} =  \k{H} (\sum_{j=1}^{N_i} -\nabla_w f_{i_j} \big / N_i)   $ from \cref{thm:compH}. We set a minibatch size of $ N_i = 256 $, which means
    that every epoch (i.e, a pass over all $N$ data pairs) contains $ 235 = \text{ceil} (N/N_i)  $ batches. For a starting vector of all zeros,
    we run stochastic gradient descent (sgd), and two compact representations with $\k{v} = \k{s} $ and  $ \k{v} = \k{y} $. Because of the stochastic 
    properties of the problem we set the initialization to be a constant identity $ \subsup{H}{k}{0} = I $ and the memory parameter to $l=1$. 
    \cref{fig:ml} displays the results of minimizing the training loss, and the accuracy for the test set.

\subsection{A Second Multiclass model}
\label{sec:ml2}
The Fashion MNIST dataset \cite{xiao2017/online} is considered to be more detailed when compared to the MNIST dataset.
Like MNIST it consists of $N=60,000$ images and labels as part of its training set and another $10,000$ pairs for a testset.
However, the images are greyscale pictures of 10 fashion items from the online retailer Zalando and can be harder to distinguish
from each other. We use a fully connected
neural network with one hidden layer (size $512 \times 512$), input layer ($ 784 \times 512 $) and output layer ($ 512 \times 10 $) 
to process each of the vectorized images with $ 784=28*28  $ pixels. We interface the compact solver with the PyTorch library 
\cite{NEURIPS2019_9015}. SGD and the compact solver with $\k{v}=\k{y}$ and $l=5$ are used to train the model. 
Since using a larger memory value can be considered as a form of regularization the compact solver can enable a larger learning rate $\alpha$.
The algorithms are tested with a minibatch size of 64, 10 epochs and a learning rate of 0.5.
\begin{table}
\centering
\setlength{\tabcolsep}{2pt}
\scriptsize
\begin{tabular}{|l l | c c c c c c c c c c|}
\hline
\multirow{2}{*}{{\small solver}} & \multirow{2}{*}{{ \small result}} & \multicolumn{10}{c | }{ { \small  epoch } } \\
& & 1 & 2 & 3 & 4 & 5 & 6 & 7 & 8 & 9 & 10 \\
\hline 
\multirow{2}{*}{ {\small  sgd }} & loss &  
0.566&	0.487&	0.466&	0.457&	0.618&	0.421&	0.406&	0.412&	0.373&	0.437 \\
& acc. & 78.40\% & 	81.50\% &	83.20\% &	84.00\% &	81.10\% &	84.60\% &	85.20\% &	85.20\% &	86.70\% &	84.90\% \\
\hline
\multirow{2}{*}{ {\small  compact }} & loss &
0.472 &	0.415 &	0.397 &	0.383 &	 0.381 &	0.365 &	0.361 &	0.368 &	0.368 &	0.367 \\
& acc. & 82.50\% &	84.60\% &	85.80\% &	86.20\% &	86.30\% &	86.90\% &	87.30\% &	87.40\% &	87.50\% &	87.80\% \\
\hline
\end{tabular}
\caption{Comparison of sgd \cite{NEURIPS2019_9015} and a compact representation algorithm on Fashion MNIST.}
\label{tab:sgdcomp}
\end{table}

\section{Conclusion}
\label{sec:concl}
This manuscript develops compact representations for two general recursive rank-2 matrix updates.
Limited-memory techniques can be efficiently implemented on top of the representations so that
computations scale linearly with the problem dimensions. By making special choices for the vector
parameters, we draw connections with known representations and enable the development of new methods
by a simple substitution for the relevant vectors. The methods scale well on large eigenvalue 
computations and appear effective for tensor factorization and regression tasks.

\appendix
\section{Proof of Corollary \ref{cor:alt}} 
\label{app:A}
\begin{proof}
In order to derive eq. \eqref{eq:compHY} in Corollary \ref{cor:alt} first represent $ \bmat{ \!\k{V} \!\!\!& \k{S} \!\!-\!\! \subsup{H}{k}{0} \k{Y}   \!} $
from eq. \eqref{eq:compthrm1}  as a product
\begin{equation*}
    \bmat{ \k{V} & \k{S} - \subsup{H}{k}{0} \k{Y}   } = \bmat{ \k{V} & \k{S} & \subsup{H}{k}{0} \k{Y}  }
    \bmat{ I    &     \\ 
                &   I \\
                & - I}
\end{equation*}
The inverse of the middle matrix in eq. \eqref{eq:compthrm1}, i.e. $ \k{M}^{-1} $, has an explicit expression as a block $ 2 \times 2 $ system
(given in eq. \eqref{eq:invM}). Note further that  the block element $ (\k{M}^{-1})_{22} = 0_{k \times k} $.  Therefore,
\begin{equation*}
    \bmat{ I    &     \\ 
                &   I \\
                & - I} 
    \bigg[ \k{M}^{-1} \bigg]
    \bmat{ I & & \\ & I & -I }
    =
    \left[
        \begin{array}{c c | c}
            & &  \\
            &  \bigg[ \k{M}^{-1} \bigg]_{\phantom{A}}  & \k{M}^{-1} \bmat { 0 \\ - I }  \\
            \hline & \bmat { 0 & - I } \k{M}^{-1}  & 0_{k \times k}
        \end{array}
        \right]
\end{equation*}
Since $ \bmat { 0 & - I } \k{M}^{-1}  $ selects the negative of the 2nd block row of $ \k{M}^{-1} $
one obtains for the product $ \bmat { 0 & - I } \k{M}^{-1}  = \bmat{ - (\subex{R}{k}{VY})^{-1} & 0 } $.
In a similar way, $ \k{M}^{-1} \bmat { 0 \\ - I }   = \bmat{  - (\subex{R}{k}{VY})^{-T} \\ 0  } $. This establishes
the formula in eq. \eqref{eq:compHY}.
\end{proof}

\section{Proof of Corollary \ref{cor:compGreen}}
\label{app:B} 
\begin{proof}
Let $ \k{V} = \k{Y}  $ in eq. \eqref{eq:compHY}, an multiple of the identity initialization 
$ H_0 = H^{(0)}_k = \gamma_k I $ and consider the block lower triangular $L$
\begin{equation*}
    L = \bmat{ I & & \\ 
                & I & \\ 
                -\frac{1}{\gamma_k}I & & I}, \quad 
    L^{-1} = \bmat{ I & & \\ 
                & I & \\ 
                \frac{1}{\gamma_k}I & & I}
\end{equation*}
Then 
\begin{equation*}
    \bmat{ \k{Y} & \k{S} & \subsc{H}{0} \k{Y}  } L = \bmat{ 0 & \k{S} & \subsc{H}{0} \k{Y} }
\end{equation*}
Moreover, $ \subex{R}{k}{VY} = \text{triu}(\kt{Y} \k{Y} ) = \subex{R}{k}{YY} $ and the expanded middle matrix 
in eq. \eqref{eq:compHY} is $\k{N} = $
\begin{align*}
    \k{N} &= \bmat{ \!\!-\!\!(\subex{R}{k}{YY})^{-T} \big( \k{R} \!\!+\!\! \kt{R} \!\!-\!\! (\k{D} \!\!+\!\! \kt{Y} \subsc{H}{0} \k{Y}  ) \big) (\subex{R}{k}{YY})^{-1}   & \!\!\! (\subex{R}{k}{YY})^{-T} & \!\!\! \!\!-\!\!(\subex{R}{k}{YY})^{\!-\!T}\! \\ 
            (\subex{R}{k}{VY})^{\!-\!1} &  0 & 0 \\
            \!\!-\!\!(\subex{R}{k}{VY})^{\!-\!1} &  0 & 0} \\
    &\equiv \bmat{ \!\!\!\!\!\!\! (\k{N})_{11} \!\!\!\! & \!\!\!\!\!\!\! (\k{N})_{21} \!\!\! & \!\!-\!(\k{N})_{21}\! \\ (\k{N})^T_{21} & 0 & 0 \\ \!\!-\!\!(\k{N})^T_{21} \!\!\!\! & 0 & 0 \\ }
\end{align*}
Then computing $ L^{-1} \k{N} L^{-T}  $ yields
\begin{equation}
    \label{eq:blockmidGreen}
    L^{-1} \k{N} L^{-T} = 
    \left[
    \begin{array}{c | c c}
    (\k{N})_{11}  &  (\k{N})_{21} & \frac{1}{\gamma_k}(\k{N})_{11} - (\k{N})_{21} \\
         \hline (\k{N})^T_{21}  & 0 &  \frac{1}{\gamma_k}(\k{N})^T_{21} \\
         \frac{1}{\gamma_k}(\k{N})_{11} - (\k{N})^T_{21} & \frac{1}{\gamma_k}(\k{N})_{21}\!\! & \frac{1}{\gamma^2_k}(\k{N})_{11} - \frac{1}{\gamma_k}(\k{N})_{21} - \frac{1}{\gamma_k}(\k{N})^T_{21}
    \end{array}
    \right]
\end{equation}
First, we consider the (3,3) block of $ L^{-1} \k{N} L^{-T}  $
\begin{align*}
    &\frac{1}{\gamma^2_k}(\k{N})_{11} \!\!-\!\! \frac{1}{\gamma_k}(\k{N})_{21} \!\!-\!\! \frac{1}{\gamma_k}(\k{N})^T_{21} = \\
    &\!-\!\frac{1}{\gamma_k}(\subex{R}{k}{YY})^{-T} \big( \k{R} \!\!+\!\! \kt{R} \!\!-\!\! (\k{D} \!\!+\!\! \kt{Y} \subsc{H}{0} \k{Y} )
    + \gamma_k\subex{R}{k}{YY} + \gamma_k\subext{R}{k}{YY} 
     \big) (\subex{R}{k}{YY})^{-1}\frac{1}{\gamma^2_k}
\end{align*}
Next we develop the the inverse of the lower $2 \times 2$ block 
\begin{align}
    \label{eq:midGreen}
    \gamma_k\bmat{
        0 &  (\k{N})^T_{21} \\
        (\k{N})_{21} & \frac{1}{\gamma_k}(\k{N})_{11} - (\k{N})_{21} - (\k{N})^T_{21}}^{\!-\!1} &= \\ 
        \bmat{  \k{R} \!\!+\!\! \kt{R} \!\!-\!\! (\k{D} \!\!+\!\! \kt{Y} \subsc{H}{0} \k{Y} )
    + \gamma_k\subex{R}{k}{YY} + \gamma_k\subext{R}{k}{YY} 
     & \!\! \gamma_k\subext{R}{k}{YY} \\ \gamma_k\subex{R}{k}{YY} & 0_{k \times k}  }
\end{align}
Using $ L L^{-1} = I $ and $ \k{V} = \k{Y} $ means that eq. \eqref{eq:compHY} becomes 
\begin{equation}
    \label{eq:compGreenLong}
    \k{H} = \subsc{H}{0} + \bmat{ 0 & \k{S} & \subsc{H}{0}\k{Y} } (L^{-1} \k{N} L^{-T}) \bmat{ 0^T \\ \kt{S} \\ \kt{Y}\subsc{H}{0} }
\end{equation}
Finally, substituting \eqref{eq:blockmidGreen} in \eqref{eq:compGreenLong}, using the inverse from \eqref{eq:midGreen} and 
relabeling the resulting middle matrix, yields the compact representation from Corollary \ref{cor:compGreen}. 
\end{proof}



\bibliographystyle{siamplain}
\bibliography{myrefs}

\begin{thebibliography}{10}

\bibitem{AcarDK11}
{\sc E.~Acar, D.~M. Dunlavy, and T.~G. Kolda}, {\em A scalable optimization
  approach for fitting canonical tensor decompositions}, Journal of
  Chemometrics, 25 (2011), pp.~67--86,
  \url{https://doi.org/https://doi.org/10.1002/cem.1335},
  \url{https://analyticalsciencejournals.onlinelibrary.wiley.com/doi/abs/10.1002/cem.1335}.

\bibitem{TensorTBX}
{\sc B.~W. Bader, T.~G. Kolda, et~al.}, {\em Tensor toolbox for matlab, version
  3.6}.
\newblock \url{www.tensortoolbox.org}, September 28, 2023.

\bibitem{BeckerLbfgs}
{\sc S.~Becker}, {\em {LBFGSB} ({L-BFGS-B}) mex wrapper}.
\newblock \url{https://www.mathworks.com/matlabcentral/
  fileexchange/35104-lbfgsb--l-bfgs-b--mex-wrapper}, 2012--2015.

\bibitem{Broyden70}
{\sc C.~G. Broyden}, {\em {The convergence of a class of double-rank
  minimization algorithms 1. General considerations}}, IMA J. Applied
  Mathematics, 6 (1970), pp.~76--90,
  \url{https://doi.org/10.1093/imamat/6.1.76},
  \url{https://doi.org/10.1093/imamat/6.1.76},
  \url{https://arxiv.org/abs/http://oup.prod.sis.lan/imamat/article-pdf/6/1/76/2233756/6-1-76.pdf}.

\bibitem{BrustBEM22}
{\sc J.~Brust, O.~Burdakov, J.~Erway, and R.~Marcia}, {\em Algorithm 1030:
  Sc-sr1: Matlab software for limited-memory sr1 trust-region methods}, ACM
  Transactions on Mathematical Software, 48 (2022), pp.~1--33.

\bibitem{Brust18}
{\sc J.~J. Brust}, {\em Large-Scale Quasi-Newton Trust-Region Methods:
  High-Accuracy Solvers, Dense Initializations, and Extensions}, PhD thesis, UC
  Merced, 2018.

\bibitem{BrustDLP21}
{\sc J.~J. Brust, Z.~Di, S.~Leyffer, and C.~G. Petra}, {\em Compact
  representations of structured bfgs matrices}, Computational Optimization and
  Applications, 80 (2021), pp.~55--88.

\bibitem{brustEM24}
{\sc J.~J. Brust, J.~B. Erway, and R.~F. Marcia}, {\em Shape-changing
  trust-region methods using multipoint symmetric secant matrices},
  Optimization Methods and Software,  (2024), pp.~1--18.

\bibitem{BrustMPS22}
{\sc J.~J. Brust, R.~F. Marcia, C.~G. Petra, and M.~A. Saunders}, {\em
  Large-scale optimization with linear equality constraints using reduced
  compact representation}, SIAM Journal on Scientific Computing, 44 (2022),
  pp.~A103--A127.

\bibitem{BurMarPil02}
{\sc O.~Burdakov, J.~Martinez, and E.~Pilotta}, {\em A limited-memory
  multipoint symmetric secant method for bound constrained optimization},
  Annals Of Operations Research, 117 (2002), pp.~51--70.

\bibitem{ByrdHNS16}
{\sc R.~H. Byrd, S.~L. Hansen, J.~Nocedal, and Y.~Singer}, {\em A stochastic
  quasi-newton method for large-scale optimization}, SIAM Journal on
  Optimization, 26 (2016), pp.~1008--1031.

\bibitem{ByrNS94}
{\sc R.~H. Byrd, J.~Nocedal, and R.~B. Schnabel}, {\em Representations of
  quasi-{N}ewton matrices and their use in limited-memory methods}, Math.
  Program., 63 (1994), pp.~129--156, \url{https://doi.org/10.1007/BF01582063}.

\bibitem{ByrdNW06}
{\sc R.~H. Byrd, J.~Nocedal, and R.~A. Waltz}, {\em Knitro: An Integrated
  Package for Nonlinear Optimization}, Springer US, Boston, MA, 2006,
  pp.~35--59, \url{https://doi.org/10.1007/0-387-30065-1_4},
  \url{https://doi.org/10.1007/0-387-30065-1_4}.

\bibitem{coleman84}
{\sc T.~F. Coleman}, ed., {\em Chapter 5 Large unconstrained optimization
  problems}, Springer Berlin Heidelberg, Berlin, Heidelberg, 1984, pp.~68--97,
  \url{https://doi.org/10.1007/3-540-12914-6_5},
  \url{https://doi.org/10.1007/3-540-12914-6_5}.

\bibitem{ConGT00a}
{\sc A.~R. Conn, N.~I.~M. Gould, and P.~L. Toint}, {\em Trust-Region Methods},
  Society for Industrial and Applied Mathematics (SIAM), Philadelphia, PA,
  2000.

\bibitem{DeGuchyEM18}
{\sc O.~DeGuchy, J.~B. Erway, and R.~F. Marcia}, {\em Compact representation of
  the full broyden class of quasi-newton updates}, Numerical Linear Algebra
  with Applications, 25 (2018), p.~e2186,
  \url{https://doi.org/https://doi.org/10.1002/nla.2186},
  \url{https://onlinelibrary.wiley.com/doi/abs/10.1002/nla.2186}.

\bibitem{deng2012mnist}
{\sc L.~Deng}, {\em The mnist database of handwritten digit images for machine
  learning research}, IEEE Signal Processing Magazine, 29 (2012), pp.~141--142.

\bibitem{DennisM77}
{\sc J.~E. Dennis, Jr and J.~J. Mor{\'e}}, {\em Quasi-newton methods,
  motivation and theory}, SIAM review, 19 (1977), pp.~46--89.

\bibitem{DennisGW81}
{\sc J.~E. Dennis~Jr, D.~M. Gay, and R.~E. Walsh}, {\em An adaptive nonlinear
  least-squares algorithm}, ACM Transactions on Mathematical Software (TOMS), 7
  (1981), pp.~348--368.

\bibitem{dennis1981adaptive}
{\sc J.~E. Dennis~Jr, D.~M. Gay, and R.~E. Walsh}, {\em An adaptive nonlinear
  least-squares algorithm}, ACM Transactions on Mathematical Software (TOMS), 7
  (1981), pp.~348--368.

\bibitem{Fletcher70}
{\sc R.~Fletcher}, {\em {A new approach to variable metric algorithms}}, The
  Computer Journal, 13 (1970), pp.~317--322,
  \url{https://doi.org/10.1093/comjnl/13.3.317},
  \url{https://doi.org/10.1093/comjnl/13.3.317},
  \url{https://arxiv.org/abs/http://oup.prod.sis.lan/comjnl/article-pdf/13/3/317/988678/130317.pdf}.

\bibitem{gillMSW84}
{\sc P.~E. Gill, W.~Murray, M.~A. Saunders, and M.~H. Wright}, {\em Sparse
  matrix methods in optimization}, SIAM Journal on Scientific and Statistical
  Computing, 5 (1984), pp.~562--589.

\bibitem{Goldfarb70}
{\sc D.~Goldfarb}, {\em {A family of variable-metric methods derived by
  variational means}}, Math. Comp., 24 (1970), pp.~23--26,
  \url{https://doi.org/10.1090/S0025-5718-1970-0258249-6},
  \url{https://doi.org/10.1090/S0025-5718-1970-0258249-6}.

\bibitem{kanzowS23}
{\sc C.~Kanzow and D.~Steck}, {\em Regularization of limited memory
  quasi-newton methods for large-scale nonconvex minimization}, Mathematical
  Programming Computation,  (2023), pp.~1--28.

\bibitem{KingmaB14}
{\sc D.~P. Kingma and J.~Ba}, {\em Adam: A method for stochastic optimization},
  arXiv preprint arXiv:1412.6980,  (2014).

\bibitem{LiuW15}
{\sc J.~Liu and S.~J. Wright}, {\em Asynchronous stochastic coordinate descent:
  Parallelism and convergence properties}, SIAM Journal on Optimization, 25
  (2015), pp.~351--376.

\bibitem{LuoEtAl15}
{\sc X.~Luo, M.~Zhou, S.~Li, Y.~Xia, Z.~You, Q.~Zhu, and H.~Leung}, {\em An
  efficient second-order approach to factorize sparse matrices in recommender
  systems}, IEEE Transactions on Industrial Informatics, 11 (2015),
  pp.~946--956, \url{https://doi.org/10.1109/TII.2015.2443723}.

\bibitem{MaBB18}
{\sc S.~Ma, R.~Bassily, and M.~Belkin}, {\em The power of interpolation:
  Understanding the effectiveness of sgd in modern over-parametrized learning},
  in International Conference on Machine Learning, PMLR, 2018, pp.~3325--3334.

\bibitem{Malouf02}
{\sc R.~Malouf}, {\em A comparison of algorithms for maximum entropy parameter
  estimation}, in Proceedings of the 6th Conference on Natural Language
  Learning - Volume 20, COLING-02, USA, 2002, Association for Computational
  Linguistics, p.~1–7, \url{https://doi.org/10.3115/1118853.1118871},
  \url{https://doi.org/10.3115/1118853.1118871}.

\bibitem{NEURIPS2019_9015}
{\sc A.~Paszke, S.~Gross, F.~Massa, A.~Lerer, J.~Bradbury, G.~Chanan,
  T.~Killeen, Z.~Lin, N.~Gimelshein, L.~Antiga, A.~Desmaison, A.~Kopf, E.~Yang,
  Z.~DeVito, M.~Raison, A.~Tejani, S.~Chilamkurthy, B.~Steiner, L.~Fang,
  J.~Bai, and S.~Chintala}, {\em Pytorch: An imperative style, high-performance
  deep learning library}, in Advances in Neural Information Processing Systems
  32, Curran Associates, Inc., 2019, pp.~8024--8035,
  \url{http://papers.neurips.cc/paper/9015-pytorch-an-imperative-style-high-performance-deep-learning-library.pdf}.

\bibitem{robbinsM51}
{\sc H.~Robbins and S.~Monro}, {\em {A Stochastic Approximation Method}}, The
  Annals of Mathematical Statistics, 22 (1951), pp.~400 -- 407,
  \url{https://doi.org/10.1214/aoms/1177729586},
  \url{https://doi.org/10.1214/aoms/1177729586}.

\bibitem{Shanno70}
{\sc D.~F. Shanno}, {\em {Conditioning of quasi-Newton methods for function
  minimization}}, Math. Comp., 24 (1970), pp.~647--656,
  \url{https://doi.org/10.1090/S0025-5718-1970-0274029-X},
  \url{https://doi.org/10.1090/S0025-5718-1970-0274029-X}.

\bibitem{MATLAB23a}
{\sc {The MathWorks Inc.}}, {\em Matlab version: 9.14.0 (r2023a)}, 2024,
  \url{https://www.mathworks.com}.

\bibitem{xiao2017/online}
{\sc H.~Xiao, K.~Rasul, and R.~Vollgraf}, {\em Fashion-mnist: a novel image
  dataset for benchmarking machine learning algorithms}, 2017,
  \url{https://arxiv.org/abs/cs.LG/1708.07747}.

\bibitem{ZhangH04}
{\sc H.~Zhang and W.~W. Hager}, {\em A nonmonotone line search technique and
  its application to unconstrained optimization}, SIAM journal on Optimization,
  14 (2004), pp.~1043--1056.

\bibitem{ZhuByrdNocedal97}
{\sc C.~Zhu, R.~Byrd, and J.~Nocedal}, {\em Algorithm 778: {L-BFGS-B}:
  {F}ortran subroutines for large-scale bound-constrained optimization}, ACM
  Trans. Math. Softw., 23 (1997), pp.~550--560.

\end{thebibliography}

\end{document}


\maketitle

\section{Proof of Corollary \cref{cor:reductbfgs}}
\label{sec:suppCorBFGS}

\begin{proof}
Let $ \k{V} = \k{S}  $ in eq. \eqref{eq:compHY}, an arbitrary symmetric initialization 
$ H_0 = H^{(0)}_k $ and consider the block lower triangular $L$
\begin{equation*}
    L = \bmat{  I   & & \\ 
                -I  & I & \\ 
                    & & I}, \quad 
    L^{-1} = \bmat{ I   & & \\ 
                    I    & I & \\ 
                         & & I}
\end{equation*}
Then 
\begin{equation*}
    \bmat{ \k{V} & \k{S} & \subsc{H}{0} \k{Y}  } L = \bmat{ \k{S} & \k{S} & \subsc{H}{0} \k{Y}  } L = \bmat{ 0 & \k{S} & \subsc{H}{0} \k{Y} }
\end{equation*}

Moreover, since $ \k{V} = \k{S} $ then  $ \subex{R}{k}{VY} =  \subex{R}{k}{SY} = \text{triu}(\kt{S} \k{Y} ) = \subsc{R}{k} $ and the expanded middle matrix 
in eq. \eqref{eq:compHY} is 
\begin{equation*}
\k{N} = 
    \bmat{ \!\!-\!\!\subsc{R}{k}^{-T} \big( \k{R} \!\!+\!\! \kt{R} \!\!-\!\! (\k{D} \!\!+\!\! \kt{Y} \subsc{H}{0} \k{Y}  ) \big) \k{R}^{-1}   & \!\!\! \k{R}^{-T} & \!\!\! \!\!-\!\!\k{R}^{\!-\!T}\! \\ 
            \k{R}^{\!-\!1} &  0 & 0 \\
            \!\!-\!\!\k{R}^{\!-\!1} &  0 & 0}
    \!\!\equiv\!\!\bmat{ \!\!\!\!\!\!\! (\k{N})_{11} \!\!\!\! & \!\!\!\!\!\!\! (\k{N})_{21} \!\!\! & \!\!-\!(\k{N})_{21}\! \\ (\k{N})^T_{21} & 0 & 0 \\ \!\!-\!\!(\k{N})^T_{21} \!\!\!\! & 0 & 0 \\ }
\end{equation*}
Then computing $ L^{-1} \k{N} L^{-T}  $ yields
\begin{equation}
    \label{eq:blockmidBFGS}
    L^{-1} \k{N} L^{-T} = 
    \left[
    \begin{array}{c | c c}
    (\k{N})_{11}  & (\k{N})_{11} + (\k{N})_{21} & - (\k{N})_{21} \\
         \hline (\k{N})_{11} + (\k{N})^T_{21} & (\k{N})_{11} + (\k{N})^T_{21} + (\k{N})_{21} &  -(\k{N})_{21} \\
         -(\k{N})^T_{21} & -(\k{N})^T_{21} & 0_{k \times k}
    \end{array}
    \right]
\end{equation}
First, using the definitions of $ (\k{N})_{11} $ and $ (\k{N})_{21} $ we consider the (2,2) block of eq. \eqref{eq:blockmidBFGS}
\begin{equation*}
    (\k{N})_{11} + (\k{N})_{21} + (\k{N})^T_{21} = 
    \k{R}^{-T} \big(\k{D} + \kt{Y} \subsc{H}{0} \k{Y} ) \k{R}^{-1}
\end{equation*}
This means that 
\begin{align}
    \k{H} &= \subsc{H}{0} + \bmat{ \k{S} & \k{S} & \subsc{H}{0} \k{Y}  } L L^{-1} \k{N} L^{-T} L^T  \bmat{ \kt{S} \\ \kt{S} \\ \kt{Y} \subsc{H}{0} } \nonumber \\
          &= \subsc{H}{0} + \bmat{ 0 & \k{S} & \subsc{H}{0} \k{Y}  } (L^{-1} \k{N} L^{-T}) \bmat{ 0^T \\ \kt{S} \\ \kt{Y} \subsc{H}{0} } \nonumber \\
          &= \subsc{H}{0} + \bmat{\k{S} & \subsc{H}{0} \k{Y}  } 
          \bmat{ (\k{N})_{11} + (\k{N})^T_{21} + (\k{N})_{21} &  -(\k{N})_{21} \\
                 -(\k{N})^T_{21} & 0_{k \times k}   } 
          \bmat{ \kt{S} \\ \kt{Y} \subsc{H}{0} } \nonumber \\
          &= \label{eq:compBFGSRED} \subsc{H}{0} + \bmat{\k{S} & \subsc{H}{0} \k{Y}  } 
          \bmat{ \k{R}^{-T} \big(\k{D} + \kt{Y} \subsc{H}{0} \k{Y} ) \k{R}^{-1} &  -\k{R}^{-T} \\
                 -\k{R}^{-1} & 0_{k \times k}   } 
          \bmat{ \kt{S} \\ \kt{Y} \subsc{H}{0} } 
 \end{align}
 The last equality in eq. \eqref{eq:compBFGSRED} is precisely the expression of the \bfgs compact representation from
 eq. \eqref{eq:compactIBFGS}, completing the proof.

\end{proof}





 






